\documentclass[reqno,11pt]{amsart}
\usepackage{amsmath, latexsym, amsfonts, amssymb, amsthm, amscd}
\usepackage[utf8]{inputenc}
\usepackage{graphics,epsf,psfrag,epsfig}
\usepackage[colorlinks=true, pdfstartview=FitV, linkcolor=blue, citecolor=blue, urlcolor=blue,pagebackref=false]{hyperref}

\usepackage[showlabels,sections,floats,textmath,displaymath]{}

\numberwithin{equation}{section}

\newtheorem{theorem}{Theorem}[section]
\newtheorem{lemma}[theorem]{Lemma}
\newtheorem{proposition}[theorem]{Proposition}
\newtheorem{cor}[theorem]{Corollary}
\newtheorem{rem}[theorem]{Remark}
\newtheorem{claim}[theorem]{Claim}

\newcommand{\ind}{\mathbf{1}}
\renewcommand{\ge}{\geq}
\renewcommand{\le}{\leq}
\newcommand{\R}{\mathbb{R}}

\newcommand{\N}{\mathbb{N}}
\renewcommand{\tilde}{\widetilde}
\renewcommand{\hat}{\widehat}

\DeclareMathSymbol{\leqslant}{\mathalpha}{AMSa}{"36} 
\DeclareMathSymbol{\geqslant}{\mathalpha}{AMSa}{"3E} 
\DeclareMathSymbol{\eset}{\mathalpha}{AMSb}{"3F}     
\renewcommand{\leq}{\;\leqslant\;}                   
\newcommand{\dd}{\,\text{\rm d}}             

\newcommand{\sumtwo}[2]{\sum_{\substack{#1 \\ #2}}} 

\newcommand{\cA}{{\ensuremath{\mathcal A}} }

\newcommand{\cN}{{\ensuremath{\mathcal N}} }

\newcommand{\cT}{{\ensuremath{\mathcal T}} }

\newcommand{\cI}{{\ensuremath{\mathcal I}} }
\newcommand{\cR}{{\ensuremath{\mathcal R}} }
\newcommand{\cV}{{\ensuremath{\mathcal V}} }

\newcommand{\bP}{{\ensuremath{\mathbf P}} }
\newcommand{\bE}{{\ensuremath{\mathbf E}} }

\newcommand{\bbE}{{\ensuremath{\mathbb E}} }

\newcommand{\bbP}{{\ensuremath{\mathbb P}} }

\newcommand{\gb}{\beta}
\newcommand{\gga}{\gamma}            
\newcommand{\gd}{\delta}
\newcommand{\gep}{\varepsilon}       
\newcommand{\gp}{\varphi}

\newcommand{\gz}{\zeta}

\newcommand{\gD}{\Delta}
\newcommand{\gk}{\kappa}
\newcommand{\go}{\omega}

\newcommand{\gO}{\Omega}
\newcommand{\gl}{\lambda}

\makeatletter
\def\captionfont@{\footnotesize}
\def\captionheadfont@{\scshape}

\long\def\@makecaption#1#2{%
  \vspace{2mm}
  \setbox\@tempboxa\vbox{\color@setgroup
    \advance\hsize-6pc\noindent
    \captionfont@\captionheadfont@#1\@xp\@ifnotempty\@xp
        {\@cdr#2\@nil}{.\captionfont@\upshape\enspace#2}%
    \unskip\kern-6pc\par
    \global\setbox\@ne\lastbox\color@endgroup}%
  \ifhbox\@ne 
    \setbox\@ne\hbox{\unhbox\@ne\unskip\unskip\unpenalty\unkern}%
  \fi
  \ifdim\wd\@tempboxa=\z@ 
    \setbox\@ne\hbox to\columnwidth{\hss\kern-6pc\box\@ne\hss}%
  \else 
    \setbox\@ne\vbox{\unvbox\@tempboxa\parskip\z@skip
        \noindent\unhbox\@ne\advance\hsize-6pc\par}%
\fi
  \ifnum\@tempcnta<64 
    \addvspace\abovecaptionskip
    \moveright 3pc\box\@ne
  \else 
    \moveright 3pc\box\@ne
    \nobreak
    \vskip\belowcaptionskip
  \fi
\relax
}
\makeatother
\def\writefig#1 #2 #3 {\rlap{\kern #1 truecm
\raise #2 truecm \hbox{#3}}}

\renewcommand{\a}{\mathrm{ann}}
\newcommand{\F}{\mathtt{F}}
\renewcommand{\P}{{\ensuremath{\mathbf P}} }
\newcommand{\E}{{\ensuremath{\mathbf E}} }
\newcommand{\Eo}{\mathbb{E}}
\newcommand{\Eot}{\widetilde{\mathbb{E}}}
\newcommand{\Pot}{\widetilde{\mathbb{P}}}
\newcommand{\Po}{\mathbb{P}}
\newcommand{\hca}{h_{c}^{\a}}
\newcommand{\Znc}{Z_{n,h_c^\a}^{\a}}

\newcommand{\Pnc}{\ensuremath{\mathbf P}_{n,h_c^{\a}}^{\a}}
\newcommand{\Zna}{Z_{n,h}^{\a}}
\renewcommand{\a}{ \mathrm{a}}
\newcommand{\Var}{{\rm Var}}
\newcommand{\la}{\langle}
\newcommand{\ra}{\rangle}
\begin{document}

\title[Hierarchical pinning model with disorder correlations]{Hierarchical pinning model in correlated random environment}
\author{Quentin Berger}


\author{Fabio Lucio Toninelli}
\address{\!\!\!\!\!\!\!Universit\'e de Lyon and CNRS\newline
Ecole Normale Sup\'erieure de Lyon,\newline
Laboratoire de Physique,\newline
46, all\'ee d'Italie,
 69364 Lyon Cedex 07 - France.\newline 
\rm {\texttt{quentin.berger@ens-lyon.fr}}\newline
\rm {\texttt{fabio-lucio.toninelli@ens-lyon.fr}}
}
\begin{abstract}
  We consider the hierarchical disordered pinning model studied in
  \cite{DHV}, which exhibits a localization/delocalization phase
  transition.  In the case where the disorder is i.i.d. (independent and
  identically distributed), the question of
  relevance/irrelevance  of disorder (i.e. whether disorder changes or not the
  critical properties with respect to the homogeneous case) is by now
  mathematically rather well understood \cite{GLT07,GLT08}. Here we
  consider the case where randomness is spatially correlated and
  correlations respect the hierarchical structure of the model; in the
  non-hierarchical model our choice would correspond to a power-law
  decay of correlations. 

In terms of the critical exponent of the homogeneous model and of the correlation decay exponent, 
we identify three regions. In the first one (non-summable correlations)
the phase transition disappears. In the second one (correlations decaying fast 
enough) the system behaves essentially like in the i.i.d. setting and 
the relevance/irrelevance criterion is not modified.
Finally, there is a region where the presence of correlations changes
the critical properties of the annealed system.

\end{abstract}

\maketitle

\section{Introduction}

A fundamental problem in the study of disordered systems is to
understand to what extent quenched (i.e. frozen) randomness modifies
the critical properties of a homogeneous (i.e. non-disordered) system.
Basically, the first question is
whether the transition survives in presence of disorder that locally
randomizes the thermodynamic parameter which measures the distance
from the critical point (e.g. for a ferromagnet $T-T_c$ can be randomized by
adding a random component to the couplings $J_{ij}$). If yes,
then one can ask whether the critical exponents are modified. The
celebrated Harris criterion \cite{Harris} states that disorder is
irrelevant (i.e. a sufficiently weak disorder does not change the
critical exponents) if $d\nu>2$, where $d$ is the space dimension and
$\nu$ is the correlation length critical exponent of the
homogeneous model, while it is relevant if $d\nu<2$. The case
$d\nu=2$ is called marginal and deciding between relevance and
irrelevance is a very model-dependent question.

Despite much effort, the Harris criterion is still far from having a
mathematical justification. In the last few years, the
\emph{disordered pinning model} \cite{GBbook,SFLN} emerged as a case
where the disorder relevance question can be attacked from a rigorous
point of view. This is a class of one-dimensional ($d=1$) models,
based on an underlying renewal process with power-law inter-arrival
distribution; the model lives in a random environment, such that the
occurrence of a renewal at step $n$ is modified with respect to the
law of the renewal by a factor $\exp(\epsilon_n)$, where $\epsilon_n$
is a sequence of i.i.d. random variables: if $\epsilon_n>0$
(resp. $\epsilon_n<0$) there is an energetic gain (resp. penalization)
in having the renewal at $n$. The pinning model exhibits a
localization/delocalization phase transition when the average $h:=\bbE
\epsilon_n$ is varied, and in the non-disordered case ($\beta^2:=Var
(\epsilon_n)=0$) the critical point $h_c$ and the critical exponent
$\nu$ can be computed exactly ($\nu$ depends only on the tail exponent
of the renewal inter-arrival law). Thanks to a series of recent works,
the Harris criterion has been put on mathematical grounds on this
case: it is now proven that, for $\beta$ small, $\nu$ does not change
if it is larger than $2$ \cite{A06,L,T} and it does change as soon as
$\beta\ne0$ if $\nu<2$ \cite{GT05}.  For the pinning model, the
relevance/irrelevance question can be also asked in the following
sense \cite{DHV}: is the critical point of the disordered model
(quenched critical point) equal to the critical point of the
\emph{annealed model}, where the partition function is replaced by its
disorder average?  It turns out that for $\beta$ small the
difference of the two critical points is zero if $\nu>2$ \cite{A06,T},
while it behaves like $\beta^{2/(2-\nu)}$ if $\nu<2$
\cite{AZ08,DGLT07}.  In the marginal case $\nu=2$, relevance of
disorder has also been shown, though in the weaker sense that the
difference between quenched and annealed critical points is non-zero
(it is essentially of order $\exp(-c/\beta^2)$, as argued in
\cite{DHV} and proven in \cite{GLT08,GLT09}).
Recently, a variational approach to the relevance/irrelevance
question, based on a large deviation principle,
has been proposed in \cite{CdH}.

Let us also add that, for the pinning model, the correlation length
exponent $\nu$ should coincide with the exponent governing the
vanishing of the free energy at the critical point: $\F(h,\beta)\simeq
(h -h_c(\beta))^\nu$ (this is proven in special situations,
e.g. \cite{G_correl,T_correl}, but it should be a rather general
fact). In the rest of this work, $\nu$ will actually denote the free
energy critical exponent.

It is widely expected, on general grounds, that correlations in the
environment may change qualitatively the Harris criterion: in the case
of a $d$-dimensional system where the correlation between the random
potentials at $i$ and $j$ decays as $|i-j|^{-\xi}$, Weinrib and Halperin  
\cite{WeinHalp83} predict that the Harris criterion is unchanged if
$\xi>d$ (summable correlations), while for $\xi<d$ the condition for
disorder irrelevance should be $\xi\nu>2$.

The study of the random pinning model with correlated disorder is
still in a rudimentary form. In \cite{Poisat} a case with finite-range
correlations was studied, and no modification of the Harris criterion
was found.  On the other extreme, in the  pinning model of
\cite{BL} not only correlations decay in a power-law way, but
potentials are so strongly correlated that in a system of length $N$
there are typically regions of size $N^b$, for some $b>0$, where the
$\epsilon_n$ take the same value. In this case, the authors of \cite{BL}
are able to compute the critical point and to give sharp estimates on
the critical behavior for $\beta>0$. In particular, they find that an
arbitrarily small amount of disorder \emph{does} change the critical
exponent, irrespective of the value of the non-disordered critical
exponent $\nu$.

Hierarchical models on diamond lattices, homogeneous or disordered
\cite{Bleher,Collet,DG}, are a powerful tool in the study of the critical behavior of
statistical mechanics models, especially because real-space
renormalization group transformations \`a la Migdal-Kadanoff are exact
in this case.  In this spirit, in the present work we consider the
hierarchical version of the pinning model introduced in the i.i.d.
setting in \cite{DHV} and later studied in \cite{GLT07,GLT08}. The
idea is to study a polymer on a diamond hierarchical lattice,
interacting with a one-dimensional defect line where the potentials
$\epsilon_n$ are placed (cf. \cite[Sec. 4.2]{DHV} and
\cite[Sec. 1.2]{GLT07} for more details on the relation with the
non-hierarchical pinning model).  Thanks to the diamond structure, the
partition function for a system of size $2^n$ turns out to be
expressed by a simple recursive relation in terms of the partition
functions of two systems of size $2^{n-1}$, cf. \eqref{eq:recur}. At
this point one can (as we will in the following) forget about the
polymer interpretation and just retain the recursion. As in the
non-hierarchical case, the system exhibits a
localization/delocalization phase transition witnessed by the vanishing
of the free energy when $h$ is smaller than a certain threshold value
$h_c(\beta)$.

We consider the case where disorder is Gaussian
and its correlation structure respects the hierarchical structure of
the model: the correlation between the potential at $i$ and $j$ is
given by $\gk^{d(i,j)}$, where $0<\gk<1$ and $d(i,j)$ is the tree
distance between $i$ and $j$ on a binary tree. The Weinrib-Halperin
criterion in this context would say that disorder is irrelevant if and
only if $\nu\log_2(1/\max(\gk, 1/2))>2$ which for $\gk=0$ (no
correlations) reduces to $\nu>2$ as for the i.i.d. case. In terms of a
parameter $B\in(1,2)$ which defines the geometry of the diamond
lattice, the criterion would read equivalently (cf. \eqref{eq:purenu})
\begin{eqnarray}
  \label{eq:WH}
\text{irrelevance}\Longleftrightarrow  \max(\gk,1/2)<B^2/4.
\end{eqnarray}

A closer inspection of the model, however, shows easily that \emph{the
  phase transition does not survive} for $\gk>1/2$ (cf. Section
\ref{sec:k12}).  When instead correlations are summable (which
corresponds to $\gk<1/2$) we find, in agreement with \eqref{eq:WH},
irrelevance if $B>\sqrt2$ (see Theorem~\ref{thm:shift} and
Proposition~\ref{prop:comparF}). As for $B\le \sqrt 2$, again we find
agreement with the Weinrib-Halperin criterion: disorder is relevant (see
Proposition~\ref{thm:smooth})
and if in addition $\gk<B^2/4$, the model behaves like in the i.i.d.  case as far as the difference between
quenched and annealed critical points is concerned, see
Theorem~\ref{thm:shift}. The crucial step (and the one which requires
the most technical work) in proving Theorem~\ref{thm:shift} (and
Proposition~\ref{prop:comparF}) is to show that for
$\gk<\min(1/2,B^2/4)$ the Gibbs measure of the annealed system near
the annealed critical point is close (in a suitable sense) to the
Gibbs measure of the homogeneous system near its critical point
(cf. Theorem~\ref{thm:boundEnc} and Proposition~\ref{prop:Pnc}).  This
requires some work, in particular because the annealed critical point
is not known explicitly for $\gk\ne0$. Once this is done, the proof of
disorder relevance/irrelevance according to $B\lessgtr \sqrt 2$ can be
obtained generalizing the ideas that were developed for the
i.i.d. model.

Finally, the region $B^2/4<\gk<1/2, B<\sqrt 2$ reserves somewhat of a
surprise: while we are not able to capture sharply the behavior of the annealed model and of the difference between quenched and annealed critical points
(as we do for $\gk<\min(1/2,B^2/4)$, see Theorem~\ref{thm:boundEnc},
Proposition~\ref{prop:Pnc} and Theorem \ref{th:diff}), we can prove that the
annealed model has
a different critical behavior than the homogeneous model with the same
value of $B$.  In particular, the contact fraction at the annealed
critical point scales qualitatively differently (as a function of the
system size) than for the homogeneous model, see
Equation~\eqref{eq:grossecorrel}.  In view of
Theorem~\ref{thm:boundEnc} mentioned above, this means that if we fix
$B<\sqrt2$ and we increase $\gk$ starting from $0$, at $\gk=B^2/4$ the
annealed system has a ``phase transition'' where its critical
properties change.  As we discuss in Section \ref{sec:k12}, this
suggests that, while for $\gk<B^2/4$ the annealed free energy near the
annealed critical point $\hca(\beta)$ has a singularity of type
$(h-\hca(\beta))^\nu$ and $\nu=\log_2/\log(2/B)$, for $B^2/4<\gk<1/2$
the annealed free energy should vanish as $h\searrow \hca(\beta)$ with
a larger exponent.

Let us conclude by discussing how our results would presumably read
for the correlated, non-hierarchical disordered pinning model.  If the
disorder is Gaussian and correlations decay as $|i-j|^{-\xi}$, then we
should get the same results as for the hierarchical model, provided
that $\log_2(1/\gk)=\xi$.  In particular, if $\xi<1$ (non-summable
correlations) there is no phase transition (the proof of
Theorem~\ref{thm:gk>1/2} can actually be easily adapted), and the
annealed system, well defined if $\xi>1$, would have a critical
behavior different from the homogeneous one if $1<\xi<2/\nu$ with
$\nu$ the free energy critical exponent of the homogeneous pinning
model.  As a side remark, let us recall that Dyson \cite{Dyson} used a
\emph{hierarchical} ferromagnetic Ising model (which, at least formally, resembles very much
our annealed pinning model, cf. \eqref{eq:discuss}) plus the
Griffiths correlation inequalities, to derive criteria for existence
of a ferromagnetic phase transition for a \emph{non-hierarchical},
one-dimensional Ising ferromagnet with  couplings decaying as
$J_{i-j}\sim |i-j|^{-\xi}$.  We stress that, in contrast, in our case
there are no available correlation inequalities which would allow to
infer directly results on the non-hierarchical pinning model starting
from the hierarchical one.


Let us now give an overview of the organization of the paper:

$\bullet$ In Section \ref{sec:model} we define the model and give preliminary results,
in particular on the homogeneous case, and we state our main results in Section \ref{sec:result};

$\bullet$ In Section \ref{sec:k12} we discuss the case $\gk>1/2$,
showing that the phase transition does not survive;

$\bullet$ In Section \ref{sec:resannealed}, we study in detail the annealed model, giving first some
preliminary tools (Section \ref{sec:prelim}), then looking at the case $\gk<1/2\wedge B^2/4$
and proving Theorem \ref{thm:boundEnc} and Proposition \ref{prop:Pnc}
 (Section \ref{sec:ann}), and finally focusing on
the case $B^2/4<\gk<1/2$ (Section~\ref{sec:parabola});

$\bullet$ In Section \ref{sec:var} we prove disorder irrelevance for $\gk<1/2,\, B>\sqrt{2}$, and
in Section \ref{sec:rel} we prove disorder relevance
for $\gk<1/2\wedge B^2/4,\, B\leq\sqrt{2}$.



\section{Model and preliminaries}
\label{sec:model}

\subsection{The hierarchical pinning model with hierarchically correlated disorder}
\label{sec:hiermodel}

Let $1<B<2$. We consider the following iteration
\begin{equation}
 Z_{n+1}^{(i)} = \frac{ Z_n^{(2i-1)}Z_n^{(2i)}+B-1}{B},
\label{eq:recur}
\end{equation}
for $n\in \N \cup \{0\}$ and $i\in\N$.
We study the case in which the initial condition is random and given by $Z_0^{(i)}=e^{\gb\go_i +h}$,
with $h\in\R$, $\gb\geq 0$ and where $\go:=\{\go_i\}_{i\in\N}$ is a sequence of centered Gaussian variables,
whose law is denoted by $\Po$.
One defines the law $\Po$ thanks to the correlations matrix $K$ and note $\gk_{ij}:=\Eo[\go_i\go_j]$.
We interpret $Z_n^{(i)}$ as the partition function
on the $i^{th}$ block of size~$2^n$.

In view of the recursive definition of the partition function, we make
the very natural choice of restricting to a correlation structure of
hierarchical type. For $p\in \mathbb N\cup\{0\}$ and $k\in\mathbb N$,
let 
\begin{equation}
  \label{eq:2}
  I_{k,p}:=\{(k-1 )2^p+1,\ldots,k 2^p\}
\end{equation}
be the $k^{th}$ block of size $2^p$. We define the hierarchical
distance $d(\cdot,\cdot)$ on $\mathbb N$ by  establishing that $d(i,j)=p$ if
$i,j$ are contained in the same block of size $2^p$  but not  in the
same block of size $2^{p-1}$. In other words, $d(i,j)$ is just the
tree distance between $i$ and $j$, if $\N$ is seen as the set of the
leaves of an infinite binary tree.

We assume that $\gk_{ij}$ depends only on $d(i,j)$ and for $d(i,j)=p$
  we write $\gk_{ij}=:\kappa_p$ with $\kappa_0=1$, $\kappa_p\ge 0$ for every $p$.
Actually, we make the explicit choice
\begin{equation}
 \gk_p=\kappa^p\quad \quad \mbox{for some} \quad \quad 0<\kappa<1/2.
\label{eq:assumpcorrel}
\end{equation} 

We will see in Section \ref{sec:k12} that the reason
why we exclude the case $\kappa\ge 1/2$
is that the model becomes less interesting (there is
no phase transition for the quenched model and the annealed model is
not well defined). For $\kappa=0$, one recovers the model with i.i.d. disorder.

It is standard that such a Gaussian law actually exists. An explicit
construction can
be obtained as follows.  Let $\mathcal I=\{I_{k,p},p\ge 0,k\in \mathbb
N\}$ and let $\{\hat\go_I\}_{I\in\cI}$ be a family
of i.i.d. standard Gaussian $\cN(0,1)$ variables, and note its law
$\hat \Po$. Then one has the following equality in law:
\begin{eqnarray}
 \label{eq:3}
  \go_i := \sum_{ I\in\cI; i\in I} \hat \gk_I \hat \go_I,
\end{eqnarray}
where $\hat
\kappa_{I_{k,p}}:=\hat\kappa_p:=\sqrt{\kappa^p-\kappa^{p+1}}$ (just
check that the Gaussian family thus constructed has the correct
correlation structure; the sum in the r.h.s. of \eqref{eq:3} is well
defined since $\sum_p\hat \kappa_p^2=1<\infty.$)

\smallskip
We point out that all our results can be easily extended to the case
where 
$\kappa:=\lim_{p\to\infty}|\gk_p|^{1/p}$
exists and is in $(0,1/2)$. 

\medskip


The \textit{quenched} free energy  of the model is defined by
\begin{equation}
\F(\gb,h):= \lim_{n\to\infty} \frac{1}{2^n} \log Z_{n,h}^{\go} 
\stackrel{\Po-a.s}= \lim_{n\to\infty} \frac{1}{2^n} \Eo[ \log Z_{n,h}^{\go} ],
\label{defF}
\end{equation}
where $Z_{n,h}^{\go}$ denotes $Z_n^{(1)}$ (it is helpful to
indicate explicitly the dependence on $h$ and on $\go$, 
the dependence on $\gb$ being implicit to get simpler notations.)
The above definition is justified by the following Theorem:
\begin{theorem}
The limit in \eqref{defF} exists $\Po$-almost surely and in $L^1(\dd \Po)$, is almost surely
constant and non-negative. The function $\F$ is convex, and $\F(\gb,\cdot)$ is non-decreasing. These properties
are inherited from
\begin{equation}
 \F_n(\gb,h):= \frac{1}{2^n} \Eo[ \log Z_{n,h}^{\go} ].
\end{equation}
 $\F_n(\gb,h)$ converges exponentially fast to $\F(\gb,h)$, and more precisely one has for all $n\geq 1$
\begin{equation}
\label{note2}
 \F_n(\gb,h) - \frac{1}{2^n} \log B \leq \F(\gb,h) \leq \F_n(\gb,h) + \frac{1}{2^n} \log \left( \frac{B^2 + B-1}{B(B-1)}\right).
\end{equation}
\label{thm:existF}
\end{theorem}

We define also the \textsl{annealed} partition
function $Z_{n,h}^{\a}:=\Eo[Z_{n,h}^{\go}]$, and the \textit{annealed} free energy:
\begin{equation}
\F^{\a}(\gb,h):= \lim_{n\to \infty} \frac{1}{2^n} \log \Eo[Z_{n,h}^{\go}].
\label{defFann}
\end{equation}

\begin{proposition}
\label{prop:existFann}
The limit in \eqref{defFann} exists, is non-negative and finite. The function $\F^{\a}$
is convex and $\F^{\a}(\gb,\cdot)$ is non-decreasing. These properties are inherited from
\begin{equation}
 \F_n^{\a}(\gb,h):= \frac{1}{2^n}  \log \Eo[ Z_{n,h}^{\go} ].
\end{equation}
$\F_n^{\a}(\gb,h)$ converges exponentially fast to $\F^{\a}(\gb,h)$, and more precisely one has for all $n\geq 1$
\begin{equation}
\label{note1}
 \F_n^{\a}(\gb,h) - \frac{1}{2^n} \log B \leq \F^{\a}(\gb,h) \leq \F_n^{\a}(\gb,h) + O((2\gk)^n).
\end{equation}
\end{proposition}
Note that the error terms in the upper bounds in \eqref{note2}-\eqref{note1} are not
of the same order.

Finiteness of the annealed free energy would fail  if the correlations
where not summable, i.e. if $\sum_j \kappa_{ij}=\infty$, which would
be the case for $\kappa\geq1/2$.

The proof of Theorem \ref{thm:existF} is almost identical to the proof
of \cite [Theorem 1.1]{GLT07} (one has just to use Kingman's
subadditive ergodic theorem instead of the law of large numbers) so we skip it.
The fact that $\F(\gb,h)<\infty$ is a trivial consequence of
$Z_{n,h}^\go\le \exp(\sum_{i=1}^{2^n}(\gb |\go_i|+h))$.
The proof of Proposition \ref{prop:existFann} is postponed to Section \ref{sec:prelim}.

We can compare the \textit{quenched} and \textit{annealed} free
energies, with the Jensen inequality:
\begin{equation}
 \F(\gb,h)= \lim_{n\to\infty} \frac{1}{2^n} \Eo[\log Z_{n,h}^{\go}] \leq 
 \lim_{n\to \infty} \frac{1}{2^n} \log \Eo[Z_{n,h}^{\go}] = \F^{\a}(\gb,h). \label{jensen free energy}
\end{equation}
The properties of $\F^{\a}$ are well known in the non-correlated case, since in this case
the annealed model is just the hierarchical
homogeneous pinning model (see the Section~\ref{sec:pure model}).
We also have the existence of
critical points for both \textsl{quenched} and \textsl{annealed} models, thanks to the convexity
and the monotonicity of the free energies with respect to $h$:

\begin{proposition}[Critical points]
 Let $\gb>0$ being fixed. There exist critical values $\hca(\gb), h_c(\gb)
$ such that
\begin{itemize}
\item $\F^{\a}(\gb,h)=0$ if $h\leq \hca(\gb)$ and $\F^{\a}(\gb,h)>0$ if $h>\hca(\gb)$

\item $\F(\gb,h)=0$ if $h\leq h_c(\gb)$ and $\F(\gb,h)>0$ if $h>h_c(\gb)$.
\end{itemize}
One  has $-c_\kappa\beta^2\le \hca(\gb)\le h_c(\gb)\le 0$ for some constant $c_\kappa<\infty$.
\end{proposition}
The inequality $\hca(\gb) \leq h_c(\gb)$ is a direct consequence of
~(\ref{jensen free energy}). The fact that $\hca(\gb)\ge -c_\kappa
\beta^2$ is discussed after \eqref{eq:discuss}.  The bound
$h_c(\gb)\le 0$ follows from $\F(\beta,h)\ge \F(0,h)$, which 
is proven in \cite[Prop. 5.1]{GBbook} (the proof is given there 
for the i.i.d. disorder model but it works identically for the 
correlated case, since it simply requires that $\Eo(\omega_i)=0$).

In the sequel, we often write $\hca$ instead of $\hca(\gb)$ for
brevity.

\subsection{Galton-Watson interpretation and polymer measure}
\label{sec:GWtree}

Let us take $1<B<2$, and set $\P_n$ the law of a Galton-Watson tree $\cT_n$ of depth $n+1$,
where the offspring distribution concentrates on $0$ with probability $\frac{B-1}{B}$
and on $2$ with probability $\frac{1}{B}$. Thus, the mean offspring size
is $2/B>1$, and the Galton-Watson process is supercritical. We then have a random binary tree with a random subset of descendants
and we define the set $\cR_n\subset  \{1,\ldots,2^n \}$ of individuals that are
present at the $n^{th}$ generation (which are the leaves of $\cT_n$).

Recall the definition \eqref{eq:2} of $I_{k,p}$, the $k^{\rm th}$ block of size $2^p$, and of
the hierarchical (tree) distance $d(\cdot,\cdot)$ introduced in Section \ref{sec:hiermodel}.

One has the useful following Proposition
\begin{proposition}[\cite{GLT08}, Proposition 4.1]
For any $n\geq 0$ and given a subset $I\subset \{1,\ldots,2^n\}$,
one defines $\cT_{I}^{(n)}$ to be the subtree of the standard binary tree of depth $n+1$,
obtained by deleting all the edges, except those which link leaves $i\in I$ to the root. We note $v(n,I)$
the number of nodes of $\cT_{I}^{(n)}$, with the convention that leaves are not counted as nodes, while the root is.
Then one has
\begin{equation}
 \E_n\left[ \gd_I \right] = B^{-v(n,I)}, 
\end{equation}
where $\gd_I:= \prod_{i\in I} \gd_i$ and
where $\gd_i=1$ if the individual $i$ is present at generation $n$
(i.e. if $i\in \cR_n$),
and $\gd_i=0$ otherwise. In particular $\E_n[\gd_i]=B^{-n}$ for every $i\in\{1,\ldots,2^n\}$.
\label{prop:nodi}
\end{proposition}

Using the recursive structure of the Galton-Watson tree $\cT_n$, one can rewrite the partition function as
\begin{equation}
 Z_{n}^{(i)} = 
\E_n \left[ \exp\left(\sum_{k=1}^{2^n} (\gb \go_{2^n(i-1)+k} +h) \gd_k \right)  \right],
\label{eq:defZ}
\end{equation}
since it satisfies the iteration \eqref{eq:recur} and the correct
initial condition $Z_0^{(i)}=\exp(\gb\go_i+h)$.
It is convenient to define 
\begin{eqnarray}
  \label{eq:4}
H_{n,h}^{\go,(i)} =
\sum_{k\in I_{i,n}} (\gb\go_k+h)\gd_k 
\end{eqnarray}
as the Hamiltonian on the
$i^{th}$ block of size $2^n$ (we also write $H_{n,h}^{\go}$ for $H_{n,h}^{\go,(1)}$ if there is no ambiguity). 
This allows to introduce the polymer measure
\begin{equation}
\frac{\dd \bP_{n,h}^{\go}}{\dd \P_n} := \frac{1}{Z_{n,h}^{\go}}
    \exp\left( H_{n,h}^{\go} \right).
\end{equation}

\begin{rem}\rm
 As in the pinning model~\cite{GBbook}, the critical point $h_c(\gb)$ marks the transition from a delocalized to a localized regime.
We observe that thanks to the convexity of the free energy, for a fixed $\gb$
\begin{equation}
 \partial_{h} \F(\gb,h) = \lim_{n\to\infty} \bE_{n,h}^{\go}\left[ \frac{1}{2^n} \sum_{k=1}^{2^n} \gd_k \right],
\end{equation}
almost surely in $\go$, for every $h$ such that $\F$ is differentiable
at $h$. This is the so-called average ``contact fraction'' under the measure
$\bP_{n,h}^{\go}$.  If $h<h_c(\beta)$, $\F(\gb,h)=0$ and the density
of contact goes to~$0$: we are in the delocalized regime. On the other
hand, if $h>h_c(\gb)$, we have $\F(\gb,h)>0$, and there is a positive
density of contacts: this is the localized regime.

Such a remark applies also naturally to the annealed model. 
\end{rem}


\subsection{Critical behavior of  the pure model}
\label{sec:pure model}

It is convenient to set 
\begin{eqnarray}
  \label{eq:5}
S_n^{(i)} = \sum_{k\in I_{i,n}}
\gd_k  
\end{eqnarray}
to be the number of contact points on the block $I_{i,n}$, and write
$S_n=S_n^{(1)}$ if there is no ambiguity. We then have of course
$S_n^{(i)}=S_{n-1}^{(2i-1)}+S_{n-1}^{(2i)}$.

The pure model is the model in which $\gb=0$: its partition function is
$Z_{n,h}^{\rm pure}=\E_n\left[ \exp\left( h S_n \right) \right]$ and we
let $\F(h)$ denote its free energy. It is well known that the pure
model exhibits a phase transition at the critical point $h_c(\gb=0)=0$:
\begin{theorem}[\cite{GLT07}, Theorem 1.2]
For every $B\in(1,2)$, there exist two constants $c_0:=c_0(B)>0$ and $c_0':=c'_0(B)>0$ such that for all
$0\leq h\leq 1$, we have 
\begin{equation}
\label{eq:pureF}
 c_0 h^{\nu} \leq \F(h) \leq c_0' h^{\nu}
\end{equation}
with
\begin{equation}
\nu = \frac{\log 2}{\log (2/B)} >1.
\label{eq:purenu}
\end{equation}
\label{thm:purebehav}
\end{theorem}
The exponent $\nu$ is called the pure critical exponent.
Note that $\nu$ is an increasing function of $B$, and that we have $\nu=2$ for $B=B_c:=\sqrt{2}$.
We give other useful estimates on the pure model in Appendix \ref{sec:purestim}.


\section{Main results}
\label{sec:result}
In this section we frequently write $h_c^a$ instead of $\hca(\gb)$.

It turns out that the effect of correlations is extremely different
according to whether 
$\gk<\frac{B^2}{4}\wedge \frac12$ or not. In the former case, our
first result says that, the correlations decaying fast enough, the critical
properties of the annealed model are very close to those of the pure
one.

First, let us write down more explicitly what $\Zna=\Eo [Z_{n,h}^{\go}]$ is.
Note that the Gaussian structure of the disorder is very helpful, to be able to give an explicit
formula for the annealed partition function, only in terms of two points correlations.
The computation gives 
\begin{equation}
 \Zna =  \E_n \left[
      \exp\left( \left(\frac{\gb^2}{2}+h\right) \sum_{k=1}^{2^n} \gd_k + 
              \gb^2/2 \sum_{p=1}^n \gk_p \sumtwo{1\leq i,j\leq 2^n}{d(i,j)=p} \gd_i \gd_j \right)
                 \right]=:\E_n \left[e^{H_{n,h}^{\a}}
\right].
 \label{eq:defZann}
\end{equation}

One easily realizes that
\begin{multline}
\label{eq:discuss}
H_{n,h}^{\a}= h  \sum_{k=1}^{2^n} \gd_k + \frac{\gb^2}{2}
\sum_{i,j=1}^{2^n} \gk_{ij}  \gd_i \gd_j  
    =   \left(\frac{\gb^2}{2}+h\right) S_n + \gb^2 \sum_{p=1}^n \gk_p \sum_{i=1}^{2^{n-p}} S_{p-1}^{(2i-1)} S_{p-1}^{(2i)}.
\end{multline}
In particular note that
\[
(h+\beta^2/2)\sum_{k=1}^{2^n}\delta_k\le H_{n,h}^{\a}\le 
(h+c_\kappa \beta^2)\sum_{k=1}^{2^n}\delta_k:=
\left(h+\frac{\beta^2}2\sum_{p\ge0}2^{p-1}\kappa^p
\right)
\sum_{k=1}^{2^n}\delta_k,
\]
which together with the fact that $h_c(\beta=0)=0$, implies
$ -c_\kappa \gb^2\le \hca(\beta)\le -\beta^2/2$.

We also use the notation $H_n^{\a,(k)}$ for the ``annealed
Hamiltonian'' on the $k^{th}$ block of size $2^n$
\[
H_{n,h}^{\a,(k)}= h  \sum_{l\in I_{k,n}} \gd_l + \frac{\gb^2}{2}
\sum_{i,j\in I_{k,n}} \gk_{ij}  \gd_i \gd_j .
\]
and the following relation holds: 
\begin{equation}
H_{n+1,h}^{\a} = H_{n,h}^{\a,(1)}+H_{n,h}^{\a,(2)}+\gb^2\gk_{n+1} S_n^{(1)}S_n^{(2)}.
\label{eq:iterHann}
\end{equation}

If we set $h=\hca+u$, so that the phase transition is at $u=0$, one has 
\begin{equation}
\Zna  =  \E_n  \left[ \exp\left( u S_n \right)
  e^{H^{\rm a}_{n,\hca}}\right]= Z_{n,\hca}^{\rm a} \E_{n,\hca}^{\rm a}  \left[ \exp\left( u S_n \right) \right],
\end{equation}
where 
\begin{equation}
\frac{\dd \P_{n,\hca}^{\rm a}}{\dd \P_n} := \frac{1}{\Znc}
    \exp\left( H_{n,\hca}^{\a} \right).
\end{equation}
The measure $\Pnc$ is the annealed polymer measure at the critical point $\hca$.

We can finally formulate our first result:
\begin{theorem}
Let $\gk<\frac{B^2}{4}\wedge \frac12$. There exist
some $\gb_0>0$ and constants $c_1,c_2>0$ such that for every $\gb\le\gb_0$ and $u\in[0,1]$, one has
\begin{equation}
-c_2\gb^2\left(\frac{4\gk}{B^2}\right)^n+ \E_n\left[ \exp\left( e^{-c_1 \gb^2} u S_n\right) \right]
   \leq \E_n\left[ \exp\left( u S_n \right) e^{H_{n,\hca}^{\rm a}} \right]
   \leq \E_n\left[ \exp\left( e^{c_1 \gb^2} u S_n \right) \right]
\end{equation}
so that, for any $u\in[0,1]$,
\begin{equation}
 \F\left(e^{-c_1 \gb^2} u\right)\leq \F^{\a}(\gb,\hca +u) \leq \F\left(e^{c_1 \gb^2} u\right).
\end{equation}
\label{thm:boundEnc}
\end{theorem}
 Theorem \ref{thm:boundEnc} is saying  that the critical behavior
of the annealed free energy around $\hca$ is the same as that of the
pure model around $h=0$ (in particular, same critical exponent $\nu$).

The essential tool is to prove that the measures $\P_n$ and
$\Pnc$ are close. This is the contents of the following Proposition:
\begin{proposition}
If $\gk<\frac{B^2}{4}\wedge \frac12$, then there exist some $\gb_0>0$ and  a constant $c_1>0$ such that, for
every $\gb\le\gb_0$, for any non-empty subset $I$ of $\{1,\ldots, 2^n\}$ one has
\begin{equation}
\left(e^{-c_1 \gb^2}\right)^{|I|}\E_n\left[\gd_I \right] \leq \E_n\left[\gd_I e^{H_{n,\hca}^{\a}}\right]
  \leq \left(e^{c_1 \gb^2}\right)^{|I|}\E_n\left[\gd_I \right],
\label{eq:uu}
\end{equation}
where $\gd_I:= \prod_{i\in I} \gd_{i}$. The case $I=\emptyset$ is
dealt with by
Lemma \ref{lem:boundZnc} below, that says that the partition function
at the critical point approaches $1$ exponentially fast: 
\begin{eqnarray}
  \label{eq:1}
e^{-c_2 \gb^2(4\gk/B^2)^n} \leq
\Znc\leq 1.  
\end{eqnarray}
\label{prop:Pnc}
\end{proposition}
Observe that \eqref{eq:1} says that if
 $\gk<\frac{B^2}{4}\wedge \frac12$ the partition function of the
annealed model at $h_c^{\a}$ is very close to that of the pure model
at its critical point $h=0$ (which equals identically~$1$).
We will see in Theorem \ref{th:parabola} that \eqref{eq:uu} fails, even for $\beta>0$ small,
if $\gk>\frac{B^2}{4}\wedge \frac12$.

With the crucial  Proposition \ref{prop:Pnc} in hand, it is not
hard to prove that 
for $\gk<\frac{B^2}{4}\wedge \frac12$  the Harris
criterion for disorder relevance is not modified by the presence of
disorder correlations:
\begin{theorem}
\label{th:diff} Let $\gk<\frac{B^2}{4}\wedge \frac12$.
\begin{itemize}
 \item If $1<B\leq B_c=\sqrt2$, then disorder is \textsl{relevant}:
the quenched and annealed critical points differ for every $\gb>0$, and:
\begin{itemize}
\item if $B<B_c$, there exist a constant $c_3>0$  such that for every $0\leq \gb\leq 1$
\begin{equation}
\label{alb}
(c_3)^{-1} \gb^{\frac{2}{2-\nu}}\leq h_c(\gb)-h_c^{\a}(\gb) \leq  c_3 \gb^{\frac{2}{2-\nu}}\ ;
\end{equation}
\item if $B=B_c$, there exist a constant $c_4>0$ and some $\beta_0>0$ 
such that for every $0\leq \gb\leq\beta_0$
\begin{equation}
\label{ulb}
\exp\left( -\frac{c_4}{\gb^4} \right)\leq h_c(\gb)-h_c^{\a}(\gb) \leq \exp\left( -\frac{c_4^{-1}}{\gb^{2/3}} \right)  .
\end{equation}
\end{itemize}

 \item If $B_c<B<2$, then disorder is \textsl{irrelevant}: there exists some $\gb_0>0$ such
that $h_c(\gb)=h_c^{\a}(\gb)$ for any $0<\gb\leq \gb_0$.
More precisely, for every $\eta>0$ and choosing $u>0$ sufficiently small,
$\F(\beta,\hca(\beta)+u)\ge (1-\eta)\F^{\rm a}(\beta,\hca(\beta)+u)$.
\end{itemize}
\label{thm:shift}
\end{theorem}
With some extra effort  one can presumably
improve the upper bound \eqref{ulb} to $\exp(-c_2^{-1}/\beta^2)$ and
the lower bound to $\exp(-c_2(\epsilon)/\beta^{2+\epsilon})$ for
every $\epsilon>0$, as is known for the uncorrelated case $\kappa=0$ 
\cite{GLT08,GLT09}. We will not pursue this line.

\begin{rem}
  \rm It is important to note that Theorems \ref{thm:boundEnc} and
\ref{thm:shift} do not
  require the knowledge of the value of $\hca$ (in general there is no
  hope to compute it exactly).  This makes the analysis of the
  quenched model considerably more challenging than in the
  i.i.d. disorder case $\kappa=0$, where it is immediate to see that
  $h_c^\a(\gb)=-\gb^2/2$.

\end{rem}

We mentioned in the introduction that for the i.i.d. model one can prove
that, when  the free-energy critical exponent $\nu$ of the homogeneous model
is smaller than $2$, such exponent is modified by an
arbitrarily small amount of disorder (more precisely, the result
is that the exponent is at least $2$ as soon as $\beta>0$).
The same holds for the model with correlated disorder:

\begin{proposition}
 If $\gk<1/2$, for every $B\in(1,2)$ there exists a constant $c(B)<\infty$ such that for all $\gb>0$ and $h\in\R$, we have
\begin{equation}
 \F(\gb,h) \leq \frac{c(B)}{\gb^2} \left( h-h_c(\gb) \right)_+ ^2.
\end{equation}
\label{thm:smooth}
\end{proposition}
We restrict to $\gk<1/2$ since otherwise there is no phase transition.

We do not give here the proof of this Proposition since, thanks to summability of the correlations, it is very
similar to the one for the i.i.d. hierarchical model \cite{LT}. 

\medskip

In the case $1/2>\gk\ge B^2/4$ correlations have a much more
dramatic effect on critical properties and in particular  we expect
them to change the
value of the annealed critical exponent from the value $\nu=\log
2/\log(2/B)$ to a larger one.
Partial results in this direction are collected in the following
Theorem, which shows that (some) critical properties of the 
annealed model differ from those of the homogenous one.

\begin{theorem}
\label{th:parabola}
  Let $B^2/4<\gk<1/2$ and $\beta>0$. In contrast with \eqref{eq:1}, the partition
  function at the critical point does not converge to $1$. Rather, one
  has 
  \begin{eqnarray}
    \label{eq:6}
    \prod_{p=0}^{n-1}Z^\a_{p,h^\a_c}\le \frac{
      1}{\beta\sqrt\gk}\left(\frac{B}{2\sqrt{\gk}} \right)^{n}.
  \end{eqnarray}
Also, the average number of individuals at generation $n$ at the
critical point satisfies
\begin{eqnarray}
  \label{eq:7}
 \E^\a_{n,h_c^\a}\left[ S_n\right]=  \E^\a_{n,h_c^\a}\left[\sum_{i=1}^{2^n}\delta_i\right]\le \frac{c(B)}\beta\frac 1{\gk^{(n+1)/2}}.
\end{eqnarray}
\end{theorem}
When proving Theorem \ref{th:parabola} we will actually
prove that the $m^{th}$ moment of $S_n$ under
$\P_{n,\hca}^\a$ is at most of
order
$\gk^{-mn/2}$. 
Therefore, with high probability $S_n$ is much smaller than
$(2/B)^n$, which would be the order of magnitude  of $S_n$ for 
$\gk<B^2/4\wedge 1/2$, as can be deduced from Propositions \ref{prop:Pnc}
and \ref{prop:nodi}.

\begin{figure}[htbp]
\centerline{
\psfrag{0}{$0$}
\psfrag{1}{$1$}
\psfrag{B}{$B$}
\psfrag{sqrt2}{$B_c=\sqrt{2}$}
\psfrag{2}{$2$}
\psfrag{12}{$\frac12$}
\psfrag{14}{$\frac14$}
\psfrag{k}{$\gk$}
\psfrag{?}{\small$\nu^{\rm pure}\stackrel{??}{<}\nu^{\a}$}
\psfrag{rel}{\small Relevant Disorder}
\psfrag{irr}{\small Irrelevant Disorder}
\psfrag{k12}{No phase transition}
\psfrag{kB24}{$\gk = \frac{B^2}{4}$}
\psfrag{nuannrel}{\small $\nu^{\rm pure}=\nu^{\a}<\nu^{\rm que}$}
\psfrag{nuannirr}{\small $\nu^{\rm pure}=\nu^{\a}=\nu^{\rm que}$}
\psfig{file=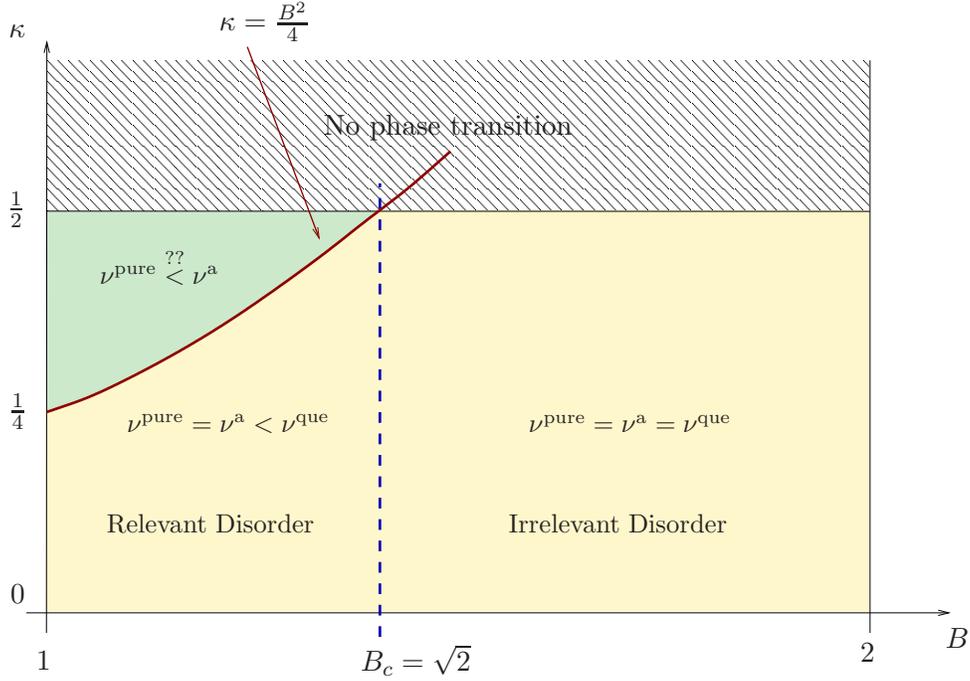,width=5in}}
  \begin{center}
    \caption{Overview of the qualitative behavior of the model. One takes $\gk<1/2$, otherwise neither annealed nor
      quenched model have any phase transition. For $\gk<1/2\wedge
      B^2/4$ the annealed model exhibits the same critical behavior as
      the pure one, and so the critical exponent is $\nu^{\a}=\nu=\log 2 \big/ \log(2/B)$. Moreover, the measures $\P_n$
      and $\P_{n,\hca}^{\rm a}$ are similar (in the sense of
      Proposition \ref{prop:Pnc}) and 
       the criterion relevance/irrelevance of disorder is the same
      as in the i.i.d. disorder case: disorder is irrelevant for
      $B>B_c:=\sqrt{2}$, marginally relevant at $B=B_c$ and relevant
      for $B<B_c$ (cf. Theorem \ref{thm:shift}).  The region above the parabola $\gk=B^2/4$ remains to
      be understood, but partial results (Theorem \ref{th:parabola})
      suggest that the critical behavior of the annealed model
      is different from the one of the pure model,  in
      particular the annealed critical exponent should be larger.
      Note that disorder is proven to be relevant for all $B<B_c$, $\gk<1/2$
      through the  ``smoothing result'' of Proposition \ref{thm:smooth},
      showing that the quenched critical exponent is strictly
      larger than the pure one.}
  \end{center}
\label{fig:zonecorrel}
\end{figure}

In other words, if we fix $B<\sqrt 2$ and we let $\gk$
grow but tuning $h$ so that we are always at the annealed critical
point, there is a phase transition in the behavior of the
finite-volume contact fraction when crossing the value $\gk=B^2/4$,
cf. also Figure \ref{fig:zonecorrel}.

\section{The case $\gk>1/2$}
\label{sec:k12}
Restricting to  the event where all the $\gd_n$ are
equal to $1$ and using Proposition \ref{prop:nodi}, one sees that
\begin{equation}
 \Zna  \geq  \left( \frac{1}{B} \right)^{2^n} \exp \left(\left((h +\gb^2/2)+\gb^2/2\sum_{p=1}^n \gk_p 2^{p-1}\right) 2^n \right).
\end{equation}
Thus, we see that 
$\F^{\a}(\gb,h)=\infty$ unless 
\begin{equation}
\label{eq:assumpK}
K_{\infty}:= \sum_{p=0}^{\infty} \gk_p 2^p < +\infty.
\end{equation}
For $\gk>1/2$, not only the annealed free energy is ill-defined. One
can also prove that the quenched free energy is strictly positive for
every value of $h\in \R$: the quenched system does not have a
localization/delocalization phase transition.
\begin{theorem}
\label{thm:gk>1/2}
If $\gk>1/2$, then $\F(\gb,h)>0$ for every $\beta>0,h\in\mathbb R$, so
that $h_c(\gb)=-\infty$. There exists some constant $c_5>0$ such that
for all $h\leq -1$ and $\gb>0$
\begin{equation}
 \F(\gb,h)\geq  \exp\left(- c_5|h|(|h|/\gb^2)^{\log 2 /\log( 2\gk)}\right).
\end{equation}
\end{theorem}
The proof of $h_c(\gb)=-\infty$ can be presumably extended to the case
$\gk=1/2$. To avoid technicalities, we do not develop this case here.

\begin{proof}
  In this proof (and in the sequel), we do not keep track of the constants $c,C,\ldots$,
  and therefore they can change from line to line.

  The idea is to lower bound the partition function by choosing a
  suitable localization strategy for the polymer to adopt, and to
  compute the contribution to the free energy of this strategy.  This
  is inspired by what is done in \cite[Chapter 6]{GBbook} to bound the
  critical point of the random copolymer model.  More precisely one
  gives a definition of a ``good block'', supposed to be favorable to
  localization in that the $\omega_i$ are sufficiently positive, and
  analyses the contribution of the strategy of aiming only at the good
  blocks.  For $\gk>1/2$ (non-summable correlations), it is a lot
  easier to find such large block (see Lemma \ref{lem:goodblock} to be
  compared with the independent case).  In this sense the behavior of
  the system is qualitatively different from the $\gk< 1/2$ case.

Clearly it is sufficient to prove the claim for $h$ negative and large enough
in absolute value.
Let us fix some $l\in\N$ (to be optimized later),
take $n>l$ and let 
$\cI\subset \{1,\dots,2^{n-l}\}$, which is supposed to denote
the set of indices corresponding to ``good blocks'' of size $2^l$.
Then for any fixed $\go$, 
targeting only the blocks in $\cI$ gives (a similar inequality was proven in
\cite{LT})
\begin{equation}
 Z_{n,h}^{\go} \geq  \left( \frac{B-1}{B^2} \right)^{v(n-l,\cI)} \prod_{k\in\cI} Z_{l,h}^{\go,(k)},
\end{equation}
where $v(n-l,\cI_n)$ is the number of nodes in the subtree
$\cT_{\cI}^{(n-l)}$ defined in Proposition \ref{prop:nodi} and
 $Z_{l,h}^{\go,(i)}$ is the partition function on $I_{k,l}$, the $i^{\rm th}$
block of size $2^l$, cf. \eqref{eq:2}.  The term $\left( \frac{B-1}{B^2} \right)^{v(n-l,\cI)}$ is
a lower bound on the probability that the node $1\le i\le 2^{n-l}$ at
generation $n-l$ has at least one descendant at level $n-l+1$ if and
only if $i\in \mathcal I$ (see Figure \ref{fig:aimgood}).
\begin{figure}[htbp]
\centerline{
\psfrag{l}{$2^l$}
\psfrag{good}{\small good blocks}
\psfig{file=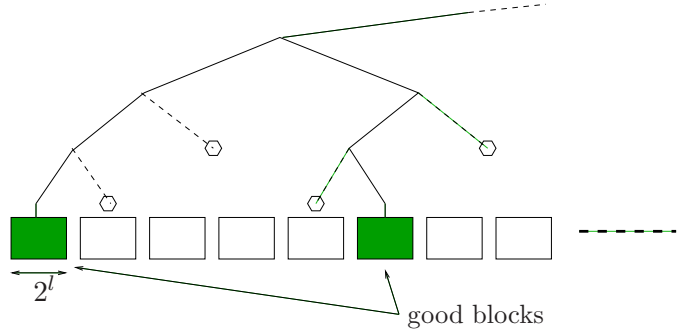,width=3.5in}}
  \begin{center}
    \caption{The strategy of aiming exactly at the good (colored)
      blocks is represented above.  One first places the subtree
      $\cT_{\cI}^{(n-l)}$, which is present with probability
      $(1/B)^{v(n-l,\cI)}$, and then forces all the leaves that do
      not lead to any good block (the hexagons in the figure) not to
      have any children, which happens with probability larger than
      $((B-1)/B)^{v(n-l,\cI)}$.  The maximal amount of nodes that
      such a tree can contain is reached when all the good blocks are
      all equally distant one from another, and is thus bounded as in
      \eqref{eq:boundnodes}.}
  \end{center}
\label{fig:aimgood}
\end{figure}

It was shown in \cite{LT} that 
\begin{equation}
 v(n,\cI) \leq |\cI| \left( 2+n-l- \lfloor\log_2 |\cI|\rfloor \right)
\label{eq:boundnodes}
\end{equation}
so that 
\begin{equation}
  \frac{1}{2^n} \log Z_{n,h}^{\go} \geq  \frac{1}{2^{n}} \sum_{k\in\cI} \log Z_{l,h}^{\go,(k)} -
  \log \left( \frac{B^2}{B-1} \right)\frac{|\cI|}{2^n} \left( 2+n-l- \lfloor\log_2 |\cI|\rfloor \right).
\label{eq:lowZhc}
\end{equation}

\medskip

Let us fix $h$ negative with $|h|$ large and take $l=l(h)\in\N$ to be chosen later.  Define then
\begin{equation}
 \cA_l^{(k)}:=\left\{ \text{for all }  i\in I_{k,l}, \text{ one has } \gb\go_i+h \geq |h|  \right\},
\end{equation}
and 
\begin{equation}
\cI(\go)=\cI_n(\go):=\{ 1\leq k\leq 2^{n-l}\ :\ \cA_l^{(k)} \text{ is verified}\}.
\end{equation}
One notices that for all $k\in\cI_n$ one has
$ Z_{l,h}^{\go,(k)} \geq  Z_{l,|h|}^{\rm pure}, $
so that one gets from \eqref{eq:lowZhc}
\begin{equation}
\frac{1}{2^n} \log Z_{n,h}^{\go} \geq   \frac{|\cI_n|}{2^{n-l}}  \frac{1}{2^l}\log Z_{l,|h|}^{\rm pure} -
     \log \left( \frac{B^2}{B-1} \right)\frac{|\cI_n|}{2^n} \left( 2+n-l- \lfloor\log_2 |\cI_n|\rfloor \right).
\end{equation}
We also note $p_l:=\Po(\cA_l^{(1)})$, so that one has
$\lim_{n\to\infty} 2^{-(n-l)} |\cI_n|=p_l$, $\bbP$-a.s., thanks to the
Ergodic Theorem.  Then, provided that $l$ is large enough so that
$2^{-l}\log Z_{l,|h|}^{\rm pure}\geq \frac12 \F(|h|)$ one has
$\bbP$-a.s.
\begin{equation}
 \F(\gb,h)\geq p_l \F (|h|)/2 -c(B) 2^{-l} p_l ( 2-\log_2 p_l )\geq p_l\left( c|h|-c' 2^{-l} (2-\log_2 p_l) \right),
\label{eq:firstlowboundF}
\end{equation}
where we used that for $|h|\geq 1$ one has $\F(|h|)\geq const\times|h|$.

\medskip
It then remains to estimate the probability $p_l$.
\begin{lemma}
\label{lemmamax}
If $\gk>1/2$, there exist two constants $c,C>0$ such that for every $l\in\N$ and $A\geq C\sqrt{l}$ one has
\begin{equation}
 \Po\left( \forall i\in\{1,\ldots,2^l\},\ \go_i\geq A \right)\geq c^{-1}\, \exp\left( -c A^2 (1/\gk)^l \right).
\end{equation}
\label{lem:goodblock}
\end{lemma}

From this lemma, and choosing $l$ such that $\sqrt{l} \leq 2|h|/(C\gb)$, one gets that
\begin{equation}
  p_l= \Po\left( \forall i\in\{1,\ldots,2^l\},\ \go_i\geq 2|h|/\gb \right)\geq c^{-1} \, \exp\left( -c \gk^{-l} h^2/\gb^2 \right).
\end{equation}
Then in view of \eqref{eq:firstlowboundF} one chooses $l= \log\left(\bar C
  |h|/\gb^2 \right)/\log(2\gk)$ (this is compatible with $\sqrt{l}
\leq 2|h|/(C\gb)$ if $|h|$ is large enough) so that 
$c|h|-c'2^{-l}(2-\log_2 p_l)\ge c|h|/2 \geq c/2$ provided that $\bar C$ is
large enough.  And \eqref{eq:firstlowboundF} finally gives with this
choice of $l$
\begin{equation}
 \F(\gb,h)\geq const\times \exp\left( - c \gk^{-l} h^2/\gb^2  \right)
    \geq const\times  \exp\left( - c' |h| \left(|h|/\gb^2\right)^{\log 2/\log(2\gk)}  \right).
\end{equation} 

\end{proof}

\begin{proof}[Proof of  Lemma \ref{lemmamax}]
First of all, note $\cA=\{\forall i\in\{1,\ldots,2^l\},\ \go_i\geq A\}$.
We consider the measure $\bar\Po$ on $\{\go_1,\ldots,\go_{2^l}\}$ which is absolutely continuous
with respect to $\Po$, and consists in translating the $\go_i$'s of $2A$, without changing the
correlation matrix $K$. Then one uses
the inequality
\begin{equation}
\Po(\cA)\geq \bar \Po(\cA) \exp\left( -\bar\Po(\cA)^{-1} ({\rm H}(\bar\Po|\Po)+e^{-1}) \right),
\label{eq:entropy}
\end{equation}
with $\rm{H}(\bar\Po|\Po)$ the relative entropy of $\bar \Po$ w.r.t. $\Po$.
Note that $\bar\Po(\cA)=\Po\left( \min\limits_{i=1,\ldots,2^l} \go_i \geq -A\right)
     = \Po\left(  \max\limits_{i=1,\ldots,2^l} \go_i \leq A \right)$,
so that from the Claim~\ref{claim:maxgaussien} below, and using that $A\geq C\sqrt{l}$, one has
$\bar\Po(\cA)\geq 1/2$.
\begin{claim}
\label{claim:maxgaussien}
Let $\{\go_i\}_{i\in\{1,\ldots,2^l\}}$ be a centered Gaussian vector
of law $\Po$, with covariance matrix $K$ such that all $\gk_{ij}\geq
0$ and $\gk_{ii}=1$. There exists a constant $C>0$ such that
\begin{equation}
 \Po\left( \max_{i=1,\ldots,2^l} \go_i \leq 
  C\sqrt{l}  \right) \geq 1/2.
\end{equation}
\end{claim}
It follows from the classical Slepian's Lemma that if
$\{\hat{\go}_i\}_{i\in\{1,\ldots,2^l\}}$ is a vector of i.i.d. standard
Gaussian variables (whose law is denoted $\hat{\Po}$), then one has
\begin{equation}
 \Eo\left[ \max_{i=1,\ldots,2^l} \go_i \right] \leq \hat{\Eo}\left[ \max_{i=1,\ldots,2^l} \hat{\go}_i \right] \leq c\sqrt{l}, 
\end{equation}
where the second inequality is classical.  Thus one gets
\begin{equation}
 \Po\left(  \max_{i=1,\ldots,2^l} \go_i \geq 2c\sqrt{l} \right) \leq
     \frac{1}{2c\sqrt{l}} \Eo\left[ \max_{i=1,\ldots,2^l} \go_i \right] \leq 1/2.
\end{equation} 

\medskip One is thus left with estimating the relative entropy
$\rm{H}(\bar\Po|\Po)$ in \eqref{eq:entropy}.  A straightforward
Gaussian computation gives \[{\rm H}(\bar\Po|\Po)= 2A^2 \la K^{-1}
\mathbf{1}, \mathbf{1}\ra\] 
where $\mathbf{1}$ is the vector whose $2^l$ elements are all equal to
$1$.  From Lemma \ref{lem:matrix} one sees that $\mathbf{1}$ is an
eigenvector of $K$, with eigenvalue $\gl:=\gk_0+\sum_{p=1}^l
2^{p-1}\gk_p\geq const\times (2\gk)^l$, so that 
${\rm H}(\bar\Po|\Po)\le  c A^2 (1/\gk)^l$, which combined with \eqref{eq:entropy} gives
the right bound.
\end{proof}

\section{Study of the annealed model}
\label{sec:resannealed}

Let us remark first of all that since $\gk_n\geq0$, thanks to \eqref{eq:iterHann} one
has $H_{n+1,h}^{\a} \geq H_{n,h}^{\a,(1)}+H_{n,h}^{\a,(2)}$, and therefore
\begin{equation}
Z_{n+1,h}^{\a}\geq \frac{(\Zna)^2+B-1}{B}.
\label{eq:recurZann}
\end{equation}
From this one deduces that $\Znc\leq 1$. Indeed, the map $x\mapsto
(x^2+(B-1))/B$ has an unstable fixed point at $1$, and $\Znc>1$ would
imply that $\F^\a(\beta,h_c^\a)>0$.








\subsection{An auxiliary partition function, proof of Proposition \ref{prop:existFann}}
\label{sec:prelim}
It is very convenient for the following to introduce a modified partition function, both for the quenched case
and for the annealed one, defining 
\begin{equation}
\bar Z_{n,h}^{\go}=\E_n \left[\exp\left(H_{n,h}^{\go}+\theta{\gb^2} \gk_n (S_n)^2
\right)\right],\indent \text{with } \theta:= \frac{\gk}{2(1-2\gk)}
\label{defZbar}
\end{equation}
and 
\begin{eqnarray}
  \label{eq:defHabar}
 \bar Z_{n,h}^{\a}=\Eo[\bar Z_{n,h}^{\go}]=\E_n \left[\exp(\bar H_{n,h}^{\a})\right],
\end{eqnarray} with
\begin{equation}
\bar H_{n,h}^{\a} = H_{n,h}^{\a} + \theta \gb^2 \gk_n
(S_n)^2.
\end{equation}
Note that $\theta$ vanishes for $\kappa\to0$ (no need of the auxiliary 
partition function for the non-correlated model) and that it diverges
for $\gk\to 1/2$, where the annealed model is not well-defined.

We also naturally define $\bar \F^{\a}(\gb,h):= \lim_{n\to\infty}
2^{-n}\log \bar Z_{n,h}^{\a}
$ (the existence of the
limit will be shown in the course of the proof of Proposition
\ref{prop:existFann}) and, using
$\delta_k\le 1$, one gets that $Z_{n,h}^{\a}\leq \bar Z_{n,h}^{\a} \leq e^{\theta\gb^2 (4\gk)^n} Z_{n,h}^{\a}$,
so that $\bar \F^{\a}(\gb,h)= \F^{\a}(\gb,h)$ (recall we chose $\gk<1/2$).
Similarly, if $\bar \F(\gb,h):= \lim_{n\to\infty}
2^{-n}\log \bar Z^\omega_{n,h}     $ then 
$\bar\F(\gb,h)=\F(\gb,h)$.

Then, from \eqref{eq:iterHann}, one gets that (recall $\gk_n=\gk^n$ and
\eqref{eq:5})
\begin{multline}
\bar H_{n+1,h}^{\a} \leq  H_{n,h}^{\a,(1)}+H_{n,h}^{\a,(2)}+\frac{\gb^2}{2}\gk^{n+1}(S_n^{(1)})^2 +\frac{\gb^2}{2}\gk^{n+1} (S_n^{(2)})^2 \\
         +  2\theta\gb^2 \gk^{n+1} (S_n^{(1)})^2+ 2\theta \gb^2 \gk^{n+1} (S_n^{(1)})^2 \\
   = H_{n,h}^{\a,(1)}+\theta \gb^2 \gk^{n} (S_n^{(1)})^2 +H_{n,h}^{\a,(2)}+\theta \gb^2 \gk^{n} (S_n^{(2)})^2 
    = \bar H_{n,h}^{\a,(1)}+\bar H_{n,h}^{\a,(2)}
\label{eq:iterHbar}
\end{multline}
where we used the self-explanatory notation $\bar H_{n,h}^{\a,(i)}$ for the
auxiliary Hamiltonian in the block $I_{i,n}$.
We used the bounds $ab\leq 1/2(a^2+b^2)$ and $(a+b)^2\leq 2(a^2+b^2)$ and then the definition of $\theta$.

This gives in particular that
\begin{equation}
\bar Z_{n+1,h}^{\a}\leq \frac{(\bar Z_{n,h}^{\a})^2+B-1}{B},
\label{eq:recurZbar}
\end{equation}
from which one deduces that $\bar Z_{n,\hca}^{\a}\geq 1$ for all $n\in\N$. Indeed,
otherwise, for some $n_0\in\N$ one has $\bar Z_{n_0,\hca}^{\a}<1$, and then one can
find some $h>\hca$ such that $\bar Z_{n_0,h}^{\a}\leq 1$, which combined with \eqref{eq:recurZbar}
gives that $\bar Z_{n,h}^{\a}\le1$ for all $n\geq n_0$. Therefore one would have
$\F^{\a}(\gb,h)=\bar \F^{\a}(\gb,h)=0$, which is  a contradiction with the definition of $\hca$.

\begin{proof}[Proof of Proposition \ref{prop:existFann}]
One has from \eqref{eq:recurZann}
\begin{equation}
 \frac{Z_{n+1,h}^{\a}}{B} \geq \left( \frac{Z_{n,h}^{\a}}{B} \right)^2,
\end{equation} 
and from \eqref{eq:recurZbar} and the fact that $\bar Z_{n,h}^{\rm a}\ge (B-1)/B$
\begin{equation}
K_B \bar Z_{n+1,h}^{\a} \leq (K_B \bar Z_{n,h}^{\a})^2 \indent \text{ with } \indent K_B=\frac{B^2+B-1}{B(B-1)}.
\end{equation} 
Therefore, the sequence $\{2^{-n} \log (Z_{n,h}^{\a}/B)\}_{n\geq 1}$
and $\{2^{-n} \log (K_B \bar Z_{n,h}^{\a})\}_{n\geq 1}$ are non-decreasing
and non-increasing respectively,
so that both converge to a limit, $\F^{\a}(\gb,h)$ and  $ \bar
\F^{\a}(\gb,h)$ respectively, but we have already remarked earlier in
this section that $\F^{\a}(\gb,h) = \bar \F^{\a}(\gb,h)$.
One finally has
\begin{equation}
\begin{split}
\label{eq:finalmente}
 \F^{\a}(\gb,h) \geq \F^{\a}_n(\gb,h) -2^{-n} \log B \\
\F^{\a}(\gb,h)= \bar \F^{\a}(\gb,h) \leq \bar \F^{\a}_n(\gb,h) +2^{-n} \log K_B, 
\end{split}
\end{equation} 
so that since $\bar \F^{\a}_n(\gb,h) \leq \F^{\a}_n(\gb,h) +
\theta\gb^2 (2\gk)^n$, one gets the desired result.

\end{proof}

\subsection{Proof of Theorem \ref{thm:boundEnc} and Proposition \ref{prop:Pnc}}
\label{sec:ann}

The really crucial point is to prove that, provided that
$\gk<\frac{B^2}{4}\wedge\frac12$, the annealed partition function (and
the auxiliary one $\bar Z_{n,h}^{\rm a}$) at the annealed critical
point converges exponentially fast to $1$.
\begin{lemma}
If $\gk<\frac{B^2}{4}\wedge\frac12$ then there exist some constant $c_2>0$ and some $\gb_0>0$
such that for any $n\geq 0$ and every $\gb\leq \gb_0$, one has
\begin{eqnarray*}
 \exp\left(-c_2 \gb^2 (4\gk/B^2)^n\right) \leq & Z_{n,\hca}^{\a} &\leq 1, \\
  1 \le& \bar Z_{n,\hca}^{\a} &\leq  \exp\left(c_2 \gb^2 (4\gk/B^2)^n\right).      
\end{eqnarray*}
\label{lem:boundZnc}
\end{lemma}

\begin{proof}[Proof of Theorem \ref{thm:boundEnc} given Lemma
  \ref{lem:boundZnc} and Proposition \ref{prop:Pnc} ]
We expand $\exp\left( u S_n \right) $, to get
\begin{equation}
 \E_n\left[ \exp\left( u S_n \right)  e^{H_{n,\hca}^{\a}}\right] 
    = \sum_{k=0}^{\infty} \frac{u^k}{k!} \E_n\left[ \left( S_n \right)^k e^{H_{n,\hca}^{\a}} \right].
\label{eq:expand}
\end{equation}
Thanks to Proposition \ref{prop:Pnc}, we have that for any $k\geq 1$
\begin{equation}
\left(e^{-c_1 \gb^2}\right)^k \E_n\left[ \left( S_n \right)^k \right]\leq \E_n\left[ \left(S_n \right)^k  e^{H_{n,\hca}^{\a}}\right]
    \leq\left(e^{c_1 \gb^2}\right)^k \E_n\left[ \left( S_n \right)^k \right],
\end{equation}
and with \eqref{eq:expand} we have then 
\begin{equation}
  \E_n\left[ \exp\left( uS_n \right) e^{H_{n,\hca}}\right] \leq
      \Znc +\E_n\left[ \sum_{k=1}^{\infty} \frac{\left(u\,e^{c_1 \gb^2}\right)^k}{k!}  \left( S_n \right)^k \right]
   \le \E_n\left[ \exp\left(e^{c_1 \gb^2}u S_n \right) \right]
\end{equation}
where we used that $\Znc\le 1$. We naturally get the other inequality in the same way
\begin{equation}
  \E_n\left[ \exp\left( uS_n \right) e^{H_{n,\hca}}\right] \geq
     \E_n\left[ \exp\left(e^{-c_1 \gb^2}u S_n \right) \right] -c_2\gb^2\left(\frac{4\gk}{B^2}\right)^n,
\end{equation}
where we used Lemma \ref{lem:boundZnc} to get that $\Znc\geq 1-c_2\gb^2(4\gk/B^2)^n $.
\end{proof}

\begin{rem}\rm
\label{rem:boundexp}
Using the same type of expansion, Proposition \ref{prop:Pnc} gives
more general results: for example, one can get
\begin{multline}
 \E_n\left[ \exp\left(e^{-pc_1 \gb^2}u (S_n)^p \right) \right]-c_2\gb^2\left(\frac{4\gk}{B^2}\right)^n\le  
       \E_n\left[ e^{H_{n,\hca}^{\a}}\exp\left( u(S_n)^p \right) \right]\\
 \E_n\left[ e^{H_{n,\hca}^{\a}}\exp\left( u(S_n)^p \right) \right] \le  \E_n\left[ \exp\left(e^{pc_1 \gb^2}u (S_n)^p \right) \right].
\end{multline}
In the sequel, we refer to this Remark to avoid repeating this kind of
computation.
\end{rem}

Before proving  Proposition \ref{prop:Pnc} and Lemma
\ref{lem:boundZnc}, we prove the following result,
valid for any $\gk<1/2$. Given $I\subset \{1,\ldots,2^n\}$ we say that
$I$ is \emph{complete} if $2i-1\in I$ for some $i\in \N$ if and only
if $2i\in I$.
\begin{lemma}
\label{presquepropPnc}
For every $n\geq 1$, and any non-empty and complete subset $I$ of $\{1,\ldots,2^n\}$, one has
\begin{eqnarray}
\left( \prod_{p=0}^{n-1} Z_{p,\hca}^{\a}\right)^{|I|} \E_n[\gd_I] & \leq  \E_n\left[ \gd_I e^{H_{n,\hca}^{\a}} \right]  & \\
     & \leq   \E_n\left[ \gd_I e^{\bar H_{n,\hca}^{\a}} \right] & 
       \leq \left( \prod_{p=0}^{n-1}\bar  Z_{p,\hca}^{\a}\right)^{|I|} \E_n[\gd_I].
\end{eqnarray} 

\end{lemma}
Note that if $I=\emptyset$, these inequalities are false, since $Z_{n,\hca}^{\a}\leq 1\leq \bar Z_{n,\hca}^{\a}$.

\begin{proof}[Proof of Lemma \ref{presquepropPnc}]
As the two bounds rely on a similar argument, that is
$H_{n+1,\hca}^{\a}\geq H_{n,\hca}^{\a,(1)}+H_{n,\hca}^{\a,(2)}$ in one case, and
$\bar H_{n+1,\hca}^{\a}\leq \bar H_{n,\hca}^{\a,(1)}+\bar H_{n,\hca}^{\a,(2)}$ in the other case,
we focus only on the lower bound.

We prove it by iteration, the case $n=1$ being trivial (the only non-empty
complete subset is $I=\{1,2\}$ and the inequalities can be checked by hand).
Now assume that the assumption is true for some $n\geq 1$ and take
$I$ a non-empty complete subset of $\{1,\ldots,2^{n+1}\}$.
We decompose $I$ into two subsets $I_1=I\cap[1,2^{n}]$ and
$I_2=I\cap[2^{n}+1,2^{n+1}]$ and we define $\tilde I_2$ to be the
subset obtained by shifting $I_2$ to the left by $2^n$. It is easy to
realize that both $I_1$ and $\tilde I_2$ are complete subsets of 
$\{1,\ldots,2^n\}$ and one has $\E_{n+1}[\gd_{I}]=
\frac{1}{B}\E_{n}[\gd_{I_1}]\E_{n}[\gd_{\tilde I_2}]$.

\begin{figure}[htbp]
\centerline{
\psfrag{i1}{$I_1$}
\psfrag{i2}{$I_2$}
\psfrag{b1}{$1^{\rm st}$ block}
\psfrag{b2}{$2^{\rm nd}$ block}
\psfrag{branch}{$1^{\rm st}$ branch = probab. $1/B$}
\psfig{file=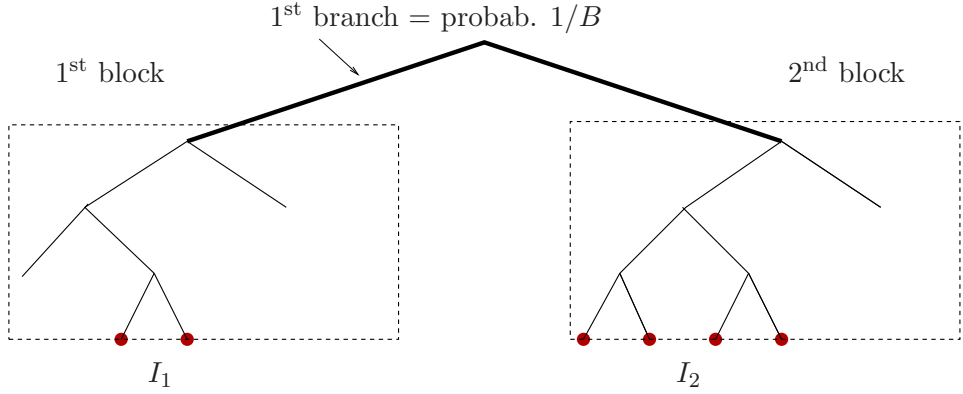,width=5in} }
  \begin{center}
    \caption{Decomposition of a non-empty complete set $I$ into two
      subsets $I_1$ and $I_2$.  If $I$ is non empty, the first
      generation must be non-empty (this has probability $1/B$). Conditionally
on this, the occurrence of $I_1$ and $I_2$ are independent events.}
  \end{center}
\label{figtree}
\end{figure}

Now, using that $H_{n+1,\hca}^{\a}\geq H_{n,\hca}^{\a,(1)}+H_{n,\hca}^{\a,(2)}$, one has
\begin{equation}
\E_{n+1}\left[\gd_{I}e^{H_{n+1,\hca}^{\a}}\right]\ge
\frac{1}{B}
\E_{n}\left[\gd_{I_1}e^{H_{n,\hca}^{\a}}\right]\E_{n}\left[\gd_{\tilde
    I_2}e^{H_{n,\hca}^{\a}}\right]
\label{iterationfacile}
\end{equation}
and two cases can occur.
\begin{enumerate}
\item 
  $\tilde I_2=\emptyset$, $|I_1|=|I|$ (or $I_1=\emptyset$, $|\tilde
  I_2|=|I|$). Then, \eqref{iterationfacile} plus the induction step gives
\begin{equation}
\E_{n+1}\left[\gd_{I}e^{H_{n+1,\hca}^{\a}}\right]\ge
\frac{1}{B} \E_{n}[\gd_{I_1}]\E_{n}[\gd_{\tilde I_2}] Z_{n,\hca}^{\a} \left( \prod_{p=0}^{n-1} Z_{p,\hca}^{\a}\right)^{|I|}.
\end{equation}
Since $Z_{n,\hca}^{\a} \leq 1$, one has $Z_{n,\hca}^{\a} \geq \left(
  Z_{n,\hca}^{\a} \right)^{|I|}$, and obtains the claim at level $n+1$.

\item  $I_1,I_2\neq \emptyset$. In this case, from \eqref{iterationfacile}, the recurrence assumption directly gives
\begin{equation}
\E_{n+1}\left[\gd_{I}e^{H_{n+1,\hca}^{\a}}\right]\ge
\frac{1}{B} \E_{n}[\gd_{I_1}]\E_{n}[\gd_{\tilde I_2}]\left( \prod_{p=0}^{n-1} Z_{p,\hca}^{\a}\right)^{|I_1|+|I_2|}.
\end{equation}
This gives the result at level $n+1$, using that $|I|=|I_1|+|I_2|$,
and bounding again $Z_{n,\hca}^{\a}\leq 1$.

\end{enumerate}
\end{proof}

\begin{proof}[Proof of Proposition \ref{prop:Pnc}]
Given $I\subset\{1,\dots,2^n\}$, let $I'$ be the smallest complete
subset of $\{1,\dots,2^n\}$ that contains $I$, and note that $|I'|\le
2|I|$. Note that 
\[
\E_n[\delta_I\exp(H^\a_{n,h_c^\a})]=\E_n[\delta_{I'}\exp(H^\a_{n,h_c^\a})],\quad\quad
\E_n[\delta_I]=\E_n[\delta_{I'}],
\]
simply because of the offspring distribution of the Galton-Watson tree: if the individual $2i-1$ is present at
generation $n$, so is the individual $2i$.
This immediately implies that the statement of Lemma
\ref{presquepropPnc} holds for every $I$ (not necessarily complete),
if $|I|$ is replaced by $2|I|$.

Then, Lemmas \ref{presquepropPnc} and \ref{lem:boundZnc} imply
Proposition \ref{prop:Pnc}
with $c_1 = 2c_2 \sum \limits_{p=0}^{\infty} \left( \frac{4\gk}{B^2} \right)^n= 2c_2 \frac{B^2}{B^2-4\gk}.$
\end{proof}

\begin{proof}[Proof of Lemma \ref{lem:boundZnc}]
  One would like to use a result analogue to Proposition
  \ref{prop:Pnc} to bound $\bar Z_{n,\hca}^{\a} = \E_n \left[
    e^{H_{n,\hca}^{\a}}\exp\left( \theta \gk^n (S_n)^2 \right)\right]
  $.  So we first prove a weaker upper bound.  The proof relies
  strongly on the pure model estimates presented in Appendix
  \ref{sec:purestim}, which show that the term $\theta \gk_n (S_n)^2$ in
  $\bar Z_{n,\hca}^{\rm a}$ has little effect if
  $\gk<\frac{B^2}{4}\wedge\frac12$.

  Take $\gp:= (2\gk) \vee \frac{4\gk}{B^2}<1$ and $C$ the constant $c$
  associated to $A=1$ in Corollary \ref{cor:bound2}, and fix some
  $\gb\leq \gb_0$, with $\gb_0:= \left(\prod_{p=0}^{\infty} e^{C
      (p+2)\varphi^p }\right)^{-2} \leq 1$.  We prove iteratively on
  $n$ that for all subsets $I$ of $\{1,\ldots, 2^n\}$ one has
\begin{equation}
\E_n \left[\gd_I e^{H_{n,\hca}^{\a}}\right] \leq  (x_n)^{|I|} \E_n[\gd_I], \indent \text{with }
 x_n := \prod_{p=0}^n e^{C (p+1)\gb \varphi^p }.
\label{firstfirstbound}
\end{equation}
Note that with our choice of $\gb_0$ one has $ (x_n)^2\leq \gb_0^{-1}$ for all $n\geq 0$.
\smallskip

The case $n=0$ is trivial (just use that $\hca\le -\beta^2/2$, as
discussed after \eqref{eq:discuss}). Now assume that \eqref{firstfirstbound} is true for some $n\geq 0$ and take
$I$ a subset of $\{1,\ldots,2^{n+1}\}$.

If $I=\emptyset$, then we simply use that $\Znc\le1$.
If $I\neq \emptyset$ decompose it as in the proof of Lemma
\ref{presquepropPnc} into two subsets $I_1,I_2$ and let $\tilde
I_2$ be obtained by translating $I_2$ to the left by $2^n$, so that
$\E_{n+1}[\gd_{I}]= \frac{1}{B} \E_{n}[\gd_{I_1}]\E_{n}[\gd_{\tilde I_2}]$ (see Figure \ref{figtree}).
Then, from the iteration \eqref{eq:iterHann} on $H_{n,h}^{\a}$ one has
\begin{equation}
H_{n+1,\hca}^{\a} \leq  H_{n,\hca}^{\a,(1)}+\frac{\gb^2}{2}\gk^{n+1} \left(S_{n}^{(1)}\right)^2 +
   H_{n,\hca}^{\a,(1)}+\frac{\gb^2}{2}\gk^{n+1} \left(S_{n}^{(2)}\right)^2,
\end{equation}
so that one gets
\begin{multline}
\E_{n+1} \left[\gd_I e^{H_{n+1,\hca}^{\a}}\right] \leq \frac{1}{B} \E_n \left[\gd_{I_1} e^{H_{n,\hca}^{\a}}
    \exp\left(\frac{\gb^2}{2}\gk^{n+1} \left(S_{n}\right)^2\right)\right]\\
     \times \E_n \left[\gd_{\tilde I_2} e^{H_{n,\hca}^{\a}} \exp\left(\frac{\gb^2}{2}\gk^{n+1} (S_{n})^2\right)\right].
\label{eq:1eredecomp}
\end{multline}

Now one can use the inductive assumption to estimate each part of
\eqref{eq:1eredecomp}. Expanding the exponential term and recalling that
$\beta_0(x_n)^2\le 1$, one has
for instance
\begin{multline}
\E_{n}\left[\gd_{I_1} e^{H_{n,\hca}^{\a}} \exp\left(\frac{\gb^2}{2}\gk^{n+1} (S_{n})^2\right)\right] =
    \sum_{k=0}^{\infty} \frac{(\gb^2 \gk^{n+1}/2)^k}{k!} \E_{n}\left[\gd_{I_1} e^{H_{n,\hca}^{\a}} \left(S_{n}\right)^{2k}\right] \\
  \leq   \sum_{k=0}^{\infty} \left(x_n\right)^{|I_1|+2k}\frac{(\gb^2 \gk^{n+1}/2)^k}{k!} \E_{n}\left[\gd_{I_1}  \left(S_{n}\right)^{2k}\right]\\
  \leq    \left(x_n\right)^{|I_1|}\E_{n}\left[\gd_{I_1}  e^{(x_n)^2\frac{\gb^2}{2}\gk^{n+1} (S_{n})^2}\right]
   \leq  \left(x_n\right)^{|I_1|}\E_{n}\left[\gd_{I_1}  e^{\frac{\gb \gk}{2}\gk^n (S_{n})^2}\right].
\label{eq:usereccur}
\end{multline}
We now use Corollary \ref{cor:bound2} to get that
\begin{equation}
  \E_{n}\left[\gd_{I_1}  e^{\frac{\gb\gk}{2}\gk^n (S_{n})^2}\right] \le
  \exp\left( C \frac{\gb\gk}{2} \gp^{-1} \gp^{n+1}\right)^{n|I_1|+1} \E_{n}[\gd_{I_1}].
\end{equation}
Combining this with \eqref{eq:1eredecomp}-\eqref{eq:usereccur} and the
definition of $\gp\geq 2\gk$ one gets
\begin{equation}
\E_{n+1} \left[\gd_I e^{H_{n+1,\hca}^{\a}}\right] \leq  (x_n)^{|I|}
     \left( e^{C \frac{\gb}{4} \gp^{n+1}} \right)^{n|I|+2} \frac{1}{B}\E_{n}[\gd_{I_1}]
\E_{n}[\gd_{\tilde I_2}].
\end{equation}
Using that  $n|I|+2\leq (n+2)|I|$ (because $I\neq \emptyset$) and the definition of $x_{n+1}=x_n e^{C(n+2)\gb \gp^{n+1} }$,
one gets equation \eqref{firstfirstbound} at level $n+1$.

\smallskip
We have performed a first crucial step: there exist some $\gb_0>0$ and a
constant $x:=\lim \limits_{n\to\infty} x_n$, such that for every $n\in\N$ and every $\gb\le\gb_0$ one has
\begin{equation}
\E_n \left[\gd_I e^{H_{n,\hca}^{\a}}\right] \leq x^{|I|} \E_n[\gd_I] \quad \text{for every } I\subset\{1,\ldots,2^n\}.
\label{step1propPnc}
\end{equation}


Then using the idea of Remark \ref{rem:boundexp},
one has from the definition of $\bar Z_{n,\hca}^{\a}$ (and expanding the exponential term)
\begin{eqnarray}
\bar Z_{n,\hca}^{\a} & = & \E_n \left[ e^{H_{n,\hca}}\right] +
\sum_{k=1}^{\infty} \frac{(\theta \gb^2 \gk^n)^k}{k!} \E_n \left[ e^{H_{n,\hca}} (S_n)^{2k}\right] \nonumber\\
   & \leq & Z_{n,\hca}^{\a}  + \E_n \left[ \exp\left(x^2 \theta \gb^2 \gk^n (S_n)^2\right)-1\right] \nonumber\\
  & \leq & Z_{n,\hca}^{\a} + \exp\left(c \gb^2 (4\gk/B^2)^n\right)-1,
\end{eqnarray}
where we used \eqref{step1propPnc} for the first inequality and Theorem \ref{thm:bound2} for the second one. 
Then using that $ Z_{n,\hca}^{\a}\leq 1$, one has the desired upper bound  for $\bar  Z_{n,\hca}^{\a}$. On the other hand,
with $\bar Z_{n,\hca}^{\a}\geq 1$ one gets that $Z_{n,\hca}^{\a}\geq 1-c' \gb^2 (4\gk/B^2)^n$, which concludes the proof.
\end{proof}

\begin{rem}\rm
\label{remarca}
  Adapting the proof of Proposition \ref{prop:Pnc} to the auxiliary
  partition function $\bar Z^\go_{n,h}$, one gets under the same
  hypothesis that there exists a constant $c_1'$ such that for any
  non-empty subset $I$ of $\{1,\ldots, 2^n\}$ one has
\begin{equation}
\left(e^{-c_1' \gb^2}\right)^{|I|}\E_n\left[\gd_I \right] \leq \E_n\left[\gd_I e^{\bar H_{n,\hca}^{\a}}\right]
  \leq \left(e^{c_1' \gb^2}\right)^{|I|}\E_n\left[\gd_I \right].
\label{eq:propPncbar}
\end{equation}
This implies, together with Lemma \ref{lem:boundZnc}, an analog of
Theorem \ref{thm:boundEnc}: there exist some $\gb_0>0$ and constants
$c_1',c_2'>0$ such that for every $\gb\le\gb_0$ and $u\in[0,1]$, one
has
\begin{equation}
\E_n\left[ \exp\left( e^{-c_1' \gb^2} u S_n\right) \right]
   \leq \E_n\left[ \exp\left( u S_n \right) e^{\bar H_{n,\hca}^{\rm a}} \right]
   \leq \E_n\left[ \exp\left( e^{c_1' \gb^2} u S_n \right) \right]+c_2'\gb^2\left(\frac{4\gk}{B^2}\right)^n.
\label{eq:thmbar}
\end{equation}
\end{rem}

\subsection{The case $B^2/4<\gk<1/2$: proof of Theorem \ref{th:parabola}}
\label{sec:parabola}

Using the identity \eqref{eq:iterHann}, one has for all $n\in\N$ and $h\in \R$
\begin{eqnarray}
Z_{n+1,h}^{\a}&=&\frac{1}{B}\E_n^{\otimes 2} \left[ e^{H_{n,h}^{\a,(1)}} e^{H_{n,h}^{\a,(2)}}
 \exp\left({\gb^2}\gk^{n+1} S_n^{(1)}S_n^{(2)}\right) \right] +\frac{B-1}{B} \\
& =& \frac{1}{B} \sum_{m=0}^{\infty} \frac{(\gb^2\gk^{n+1})^m}{m!}
  \E_n\left[e^{H_{n,h}^{\a}} (S_n)^m\right]^2 + \frac{B-1}{B}.
\label{iterannexact}
\end{eqnarray}
If one takes $h=\hca$ and uses the bound $Z_{n+1,h_c^\a}^{\a}\le 1$, one gets
\begin{equation}
 \sum_{m=0}^{\infty} \frac{(\gb^2\gk^{n+1})^m}{m!}
  \E_n\left[e^{H_{n,\hca}^{\a}} (S_n)^m\right]^2\leq 1,
\end{equation}
so that bounding each term of the sum by $1$, one gets that for all $m\geq 0$
\begin{equation}
  \E_n\left[e^{H_{n,\hca}^{\a}} (S_n)^m\right]\leq \sqrt{m!}\left(\frac{1}{\gb} \left(\frac{1}{\sqrt{\gk}}\right)^{n+1}\right)^m.
\label{eq:grossecorrel}
\end{equation}
For $m=1$ (using $Z_{n,h}^{\rm a}\ge (B-1)/B$) we obtain
\eqref{eq:7}, but  also  an estimate for all
the moments of $S_n$.

Using Lemma \ref{presquepropPnc} one has
\begin{equation}
\left(\frac{2}{B}\right)^{n}\prod_{p=0}^{n-1} Z_{p,\hca}^{\a} \leq \E_n\left[e^{H_{n,\hca}^{\a}} S_n\right] 
     \leq \frac{1}{\gb} \left(\frac{1}{\sqrt{\gk}}\right)^{n+1},
\end{equation}
which implies \eqref{eq:6}.
Another observation is that, writing $h=\hca+u$, one gets from \eqref{eq:grossecorrel} that
\begin{equation}
 \E_n\left[ e^{H_{n,h}^{\a}}\right] =  \E_n\left[e^{uS_n} e^{H_{n,\hca}^{\a}}\right]
  \leq \sum_{m=0}^{\infty} \frac{1}{\sqrt{m!}}  \left(\frac{u}{\gb}  \left(\frac{1}{\sqrt{\gk}}\right)^{n+1}\right)^m.
\end{equation} 
Thus if $u\leq \left(\sqrt{\gk}\right)^{n}$, one has that $
Z^\a_{n,\hca+u}=\E_n\left[ e^{H_{n,h}^{\a}}\right]$ does not grow with $n$.
This is in contrast with the pure model 
where
\[
Z_{n,u}^{\rm pure}=\E_n[\exp(u S_n)]\ge \exp(u\E_n(S_n))=\exp(u (2/B)^n)
\]
which diverges with $n$ if $u=(\sqrt\gk)^n$ (recall we are considering 
$\gk>B^2/4$).

All these facts lead us to conjecture that the phase transition of the
annealed model
for $B^2/4<\gk<1/2$ is smoother than that of the pure model.

\section{Disorder irrelevance}
\label{sec:var}

To prove disorder irrelevance for $B>B_c$ and the upper bounds 
on the difference between quenched and annealed critical points in Theorem
\ref{thm:shift}, we use the following Proposition:
\begin{proposition}
Let $\gk<(B^2/4\wedge 1/2)$.
If $B>B_c$, there exists a $\gb_0>0$ such that
for $\gb\leq \gb_0$ and
for every $\eta \in (0,1)$ one can find
$\gep>0$ such that for all $u\in(0,\gep)$
\begin{equation}
 \F(\gb,h_c^{\a}+u) \geq (1-\eta)\F^{\a}(\gb,h_c^{\a}+u).
\end{equation}

If $B<B_c$, then
for every $\eta \in (0,1)$ one can find
 constants $c,\beta_0,\epsilon> 0$ such that if $\gb\leq \beta_0$,
for all $u\in(c \gb^{\frac{2}{2-\nu}},\epsilon(\eta) )$
\begin{equation}
 \F(\gb,h_c^{\a}+u) \geq (1-\eta)\F^{\a}(\gb,h_c^{\a}+u)
\end{equation}
with $\nu$ as in \eqref{eq:purenu}. 

If $B=B_c$, then for every  $\eta\in(0,1)$
one can find  $\beta_0>0$ and a constant $c>0$ such that if $\gb\leq \beta_0$,
for all $u\in(c \exp\left( -c \gb^{-2/3} \right), 1)$
\begin{equation}
 \F(\gb,h_c^{\a}+u) \geq (1-\eta)\F^{\a}(\gb,h_c^{\a}+u).
\end{equation}
\label{prop:comparF}
\end{proposition}


\begin{proof}
 
This is based on the study of the variance $\cV_n:=\Eo[(\bar Z_{n,h}^{\go})^2]-\Eo[\bar Z_{n,h}^{\go}]^2$.

Fix some $B\in(1,2)$.
One has
\begin{eqnarray}
 \Eo\left[ \left(\bar Z_{n,h}^{\go} \right)^2\right] =
     \E_n^{\otimes 2} \left[ \exp\left( \bar H_{n,h}^{\a}(\gd)+\bar H_{n,h}^{\a}( \gd') +
\gb^2\sum_{i,j=1}^{2^n} \gk_{ij}\gd_i\gd'_j\right)\right]
\end{eqnarray}
with $\gd$ and $\gd'$ two independent copies of the same Galton-Watson process.
We also have $ \Eo\left[\bar Z_{n,h}^{\go}\right]^2 = 
   \E_n^{\otimes 2} \left[ \exp\left(\bar H_{n,h}^{\a}(\gd)+ \bar H_{n,h}^{\a}( \gd')\right)\right]$.
To simplify notations, we write $h=\hca+u$ and we define
\begin{equation}
\label{defDn}
 D_n : = \sum_{i,j=1}^{2^n} \gk_{ij}\gd_i\gd'_j.
\end{equation} 
Then,
\begin{multline}
 \cV_n=\E_n^{\otimes 2} \left[ e^{uS_n} e^{u S'_n}\left(e^{\gb^2 D_n}-1\right)
       e^{\bar H_{n,\hca}^{\a}(\gd)}e^{\bar H_{n,\hca}^{\a}(\gd')}\right] \\
   \leq \tilde \cV_n := \E_n^{\otimes 2} \left[ e^{CuS_n} e^{CuS'_n}
        \left(e^{C\gb^2 D_n}-1\right)\right],
\label{variancebound1}
\end{multline}
where  we expanded the exponential and used Remark \ref{rem:boundexp}
and Eq. \eqref{eq:propPncbar}.

Using the Cauchy-Schwarz inequality in \eqref{variancebound1}, 
\begin{equation}
 \tilde \cV_n \leq \E_n \left[ e^{2CuS_n}\right]
  \sqrt{ \E_n^{\otimes 2} \left[\left(e^{C\gb^2 D_n}-1\right)^2\right]}
   \leq \E_n \left[ e^{2CuS_n}\right]
        \sqrt{\E_n^{\otimes 2} \left[ e^{2C\gb^2 D_n}-1\right]}.
\end{equation}
We define  $Q_n:=\cV_n/\Eo[\bar Z_{n,h}^{\go}]^2\leq  \cV_n$,
(recall that  $h\geq \hca$ and that $\Eo\bar Z_{n,\hca}^{\go}\geq 1$).  Then one also uses
Lemma \ref{lem:estimEsum} to get that $\E_n \left[
  e^{2CuS_n}\right]\leq c \exp\left( c2^n u^{\nu}
\right)$. Therefore, one has 
\begin{equation}
 Q_n \leq c  \exp\left( c2^n u^{\nu} \right)\sqrt{\E_n^{\otimes 2} \left[ e^{2C\gb^2 D_n}-1\right]}.
\end{equation} 
Defining 
\begin{equation}
\label{eq:n1}
  n_1 = n_1(u):= \log (1/u) /\log(2/B) = \nu\log (1/u)/\log 2 ,
\end{equation}
which is the value of $n$ at which $\E_n[\exp(uS_n)]$ starts
getting large,
one has for  $p\geq 0$
\begin{equation}
 Q_{n_1+p}\leq ce^{c2^{p}} \sqrt{\E_{n_1+p}^{\otimes 2} \left[ e^{2C\gb^2 D_{n_1+p}}-1\right]}.
\label{boundQ}
\end{equation} 
Thus it is left to estimate the last term, with Proposition \ref{prop:bound22}.

\subsection{The case $B>B_c$}
Thanks to Proposition \ref{prop:bound22} there exists some $\gb_0>0$ such that
for  $\gb<\gb_0$ and for all $n\in\N$
\begin{equation}
 \E_n^{\otimes 2}\left[ e^{2C\gb^2 D_n}-1\right] 
 \leq c\gb_0^2 \Phi^n,
\end{equation}
for some $\Phi<1$.
Choose $p_1=p_1(n_1)$ such that 
$e^{c 2^{p_1}}\sqrt{\Phi^{n_1}}=1$ 
(note that $p_1$ diverges with $n_1$) and then
\begin{equation}
 Q_{n_1+p_1}\leq c' \sqrt{{\Phi}^{ p_1}} \stackrel{n_1\to\infty}{\longrightarrow} 0.
\label{Qok1}
\end{equation} 

Then we use that
\begin{equation}
\label{eq:thenwe}
 \Eo\left[\log \bar Z^\go_{n,h}\right] \geq \log \left( \frac{\Eo[\bar Z^\go_{n,h}]}{2}\right) 
     \Po\left( \bar Z^\go_{n,h} \geq \frac{\Eo[\bar Z^\go_{n,h}]}{2}\right)
    + \log\left( \frac{B-1}{B}\right),
\end{equation}
where $\Po\left(\bar Z^\go_{n,h} \geq \Eo[\bar Z^\go_{n,h}]/2\right) \geq 1-4Q_n$ from the
Tchebyshev inequality.
We apply this with $n=n_1+p_1(n_1)$ to get
(using also Theorem \ref{thm:existF} and \eqref{eq:finalmente})
\begin{multline}
 \F(\gb,h)\geq  \frac{1}{2^{n}}  \Eo\left[\log \bar Z^\go_{n,h}\right]- \frac{\log B}{2^{n}} \geq  
          (1-4\eta) \frac{1}{2^{n}}\log \left( \Eo[\bar Z^\go_{n,h}]\right) -\frac{c}{2^{n}} \\
     \geq  (1-4\eta) \F^{\a}(\gb,h) - \frac{c'}{2^{p_1(n_1)}} 2^{-n_1} \geq (1-5\eta)\F^{\a}(\gb,h),
\label{lowFque}
\end{multline}
provided that $n_1$ is large enough to ensure both
\begin{eqnarray}
                        Q_{n_1+p_1}&\leq& c' {\Phi}^{ p_1 (n_1)/2}\leq \eta \\
\text{and}\quad c'2^{-p_1(n_1)} u^{\nu}&\leq &\eta\F^{\a}(\gb,h) \hspace{1.5cm}\mbox{ for all}\quad u\in(0,1). \label{eq:p1u}
\end{eqnarray}
Note that the requirement on $n_1$ in \eqref{eq:p1u} also depends only on $\eta$,
cf. Theorem \ref{thm:boundEnc}. Since $n_1$ is related to $u$ via 
\eqref{eq:n1}, one has actually to assume that $u\le \epsilon(\eta)$
with $\epsilon$ sufficiently small, as required in Proposition 
\ref{prop:comparF}.

\subsection{The case $B< B_c$}
Given $\eta>0$ and $\beta\le1$, fix some $p_1=p_1(\eta)$ 
such that \eqref{eq:p1u} holds
and assume that
$c_1\beta^{2/(2-\nu)}\le u\le \epsilon(\eta)$ with $c_1=c_1(\eta)$ to be chosen
sufficiently large later (observe that if $\epsilon(\eta)$ is small one has 
that $n_1 $ and $p_1$ are large, so the above requirement 
on $p_1$ is coherent).
The definition of $n_1(u)$ (which gives $u=(B/2)^{n_1}$) and of $\nu$ (which gives $(2/B)^{\nu}=2$)
imply that 
\begin{eqnarray}
  \label{eq:n1beta}
  \beta^2\le c_1^{-1}\left(\frac2{B^2}\right)^{p_1(\eta)} \left(\frac{B^2}2\right)^{n_1+p_1(\eta)}
       \leq c_2\left(\frac{B^2}2\right)^{n_1+p_1(\eta)}
\end{eqnarray}
where $c_2=c_2(\eta)$ can be made arbitrarily small by choosing $c_1$
large.  Then, again provided that $c_2$ is small enough (i.e.  $c_1$
large enough), we can apply Proposition \ref{prop:bound22} to get from
\eqref{boundQ}
\begin{eqnarray}
  Q_{n_1+p_1(\eta)}\le c\,e^{c 2^{p_1(\eta)}} \sqrt{c \gb^2 \left(\frac2{B^2}\right)^{n_1+p_1(\eta)}}
    \le c'\,e^{c 2^{p_1(\eta)}} \sqrt{c_2(\eta)}\leq \eta.
\end{eqnarray}
From this point on, the proof proceeds like in the case $B>B_c$, 
starting from \eqref{eq:thenwe}.

\subsection{The case $B=B_c$}
This is similar to the case $B<B_c$. The value of $\beta_0$ has to be chosen
small enough to guarantee that Proposition \ref{prop:bound22} is
applicable. We skip details.
 \end{proof}

\section{Disorder relevance: critical point shift lower bounds}
\label{sec:rel}

To prove disorder relevance, we give a finite size condition for
delocalization, adapting the fractional moment method, first used
in \cite{DGLT07}, and then in \cite{GLT08,GLT09} for the pinning model
with i.i.d. disorder.

\subsection{Fractional moment iteration}
\label{sec:fracmom}


For $\gamma<1$  let $x_{\gga}$ to be the largest solution of
$$x= \frac{ x^2+(B-1)^{\gga}}{B^{\gga}}.$$
One can easily see that for $\gamma$ sufficiently close to $1$ (which we assume to be the case in the following) $x_{\gga}$ actually
exists and is strictly less
than $1$. Moreover one has that $x_\gamma$ increases to $1$ as $\gamma$ increases to $1$.
Then we have:
\begin{proposition}
Take $\gk<1/2$.
Then, setting $A_n:=\Eo\left[ \left(\bar Z^\go_{n,h}\right)^{\gga}\right]$ 
with $\bar Z^\go_{n,h}$ defined in \eqref{defZbar}, one has
\begin{equation}
 A_{n+1}\leq \frac{ A_n^2+(B-1)^{\gga}}{B^{\gga}}.
\label{itermomfrac}
\end{equation} 
If there exists some $n_0$ such that $A_{n_0} \leq x_{\gga}$, then $\F(\gb,h)=0$.
\end{proposition}

\begin{proof}
%
  If for some $n_0$ one has $A_{n_0}\leq x_{\gga}$, then iterating
  \eqref{itermomfrac} one gets $A_{n}\leq x_{\gga}\leq 1$ for all
  $n\geq n_0$.  Using the Jensen's inequality one has
\begin{equation}
\frac{1}{n} \Eo[\log \bar Z^\go_{n,h}] = \frac{1}{\gga n} \Eo[\log (\bar Z^\go_{n,h})^{\gga}]
        \leq \frac{1}{\gga n} \log A_{n}
\end{equation}
which gives $\F(\gb,h)=\bar \F(\gb,h)=0$ (equality of the two
free energies was noted after \eqref{defZbar}).

We now turn to the proof of \eqref{itermomfrac}. We define
$Z_{n,h}^{\mu}=\E_n\left[e^{H_{n,h}^\go} e^{\mu \gk_n \gb^2 (S_n)^2} \right] $
and use that $(S_{n+1})^2\leq 2(S_n)^2+ 2(S_n)^2$ to get the iteration
\begin{equation}
  Z_{n+1,h}^{\mu} \leq \frac{1}{B} Z_{n,h}^{2\gk\mu,(1)}Z_{n,h}^{2\gk\mu,(2)} + \frac{B-1}{B}
\end{equation} 
where as usual the two partition functions in the r.h.s. refer to the
first and second sub-system of size $2^n$.  From this, and using the
inequality $(a+b)^{\gga}\leq a^{\gga}+b^{\gga}$ for any $a,b\geq0$ and
$\gga \leq 1$, one has
\begin{equation}
 \Eo\left[(Z_{n+1,h}^{\mu})^{\gga} \right] \leq \frac{1}{B^{\gga}}
  \Eo\left[ \left(Z_{n,h}^{2\gk\mu ,(1)}Z_{n,h}^{2\gk\mu,(2)}\right)^{\gga}\right]  + \frac{(B-1)^{\gga}}{B^{\gga}}.
\end{equation} 
One then shows the following
\begin{lemma}
\label{lem:interpolon}
 If $\mu\geq \theta $ with $\theta = \frac{\gk}{2(1-2\gk)}$ as
in \eqref{defZbar}, 
\begin{equation}
\Eo\left[ \left(Z_{n,h}^{2\mu\gk ,(1)}Z_{n,h}^{2\mu \gk,(2)}\right)^{\gga}\right] \leq \Eo\left[\left(Z_{n,h}^{\mu} \right)^{\gga} \right] ^2.
\end{equation} 
\end{lemma}
This  gives directly \eqref{itermomfrac}, taking $\mu=\theta $ so that $ Z_{n,h}^{\mu}=\bar Z^\go_{n,h}$.
\end{proof}

\begin{proof}[Proof of Lemma \ref{lem:interpolon}]
One sets
\begin{equation}
\Phi(t,\mu):= \log \Eo_t\left[ \left(Z_{n,h}^{\mu ,(1)}Z_{n,h}^{\mu,(2)}\right)^{\gga} \right],
\end{equation}
where one defines $\Po_t$ to be the law of a Gaussian vector $(\go_1,\ldots,\go_{2^{n+1}})$ with correlations 
$\gk_{ij}(t)=\gk_p$ if $d(i,j)=p\le n$, and $\gk_{ij}(t)=t\gk_{n+1}$ if $d(i,j)=n+1$.
Then one can compute  the derivatives of $\Phi$. Using the definition of $Z_{n,h}^{\mu}$ one has for  $t\geq 0$, $\mu\in\R$
\begin{multline}
 \frac{\partial \Phi}{ \partial \mu} (t,\mu) =
      \frac{\gamma\gk_n \gb^2}{ \Eo_t\left[ \left(Z_{n,h}^{\mu ,(1)}Z_{n,h}^{\mu,(2)}\right)^{\gga}\right]} \\
  \times \Eo_t \left[\E_n^{\otimes 2} \left[ \left( (S_n^{(1)})^2 + (S_n^{(2)})^2 \right)
            e^{H_{n,h}^{\go,(1)}+H_{n,h}^{\go,(2)}} e^{\mu \gk_n ((S_n^{(1)})^2+(S_n^{(2)})^2)}\right] 
     \left(Z_{n,h}^{\mu ,(1)}Z_{n,h}^{\mu,(2)}\right)^{\gga-1}\right].
\end{multline} 
Thanks to Proposition \ref{prop:derivEot} one gets
\begin{equation}
 \frac{\partial \Phi}{ \partial t} (t,\mu) = \frac{\gk_{n+1}}{ \Eo_t\left[ \left(Z_{n,h}^{\mu ,(1)}Z_{n,h}^{\mu,(2)}\right)^{\gga} \right]}
   \sum_{i=1}^{2^{n}}\sum_{j=2^n+1}^{2^{n+1}}
     \Eo_t \left[\frac{\partial^2}{\partial \go_i \partial\go_j} \left(Z_{n,h}^{\mu ,(1)}Z_{n,h}^{\mu,(2)}\right)^{\gga}\right].
\end{equation} 
For the values of $i,j$ under consideration one has

\begin{multline}
 \frac{\partial}{\partial \go_i \partial\go_j} \left(Z_{n,h}^{\mu ,(1)}Z_{n,h}^{\mu,(2)}\right)^{\gga}
 = \gga ^2 \gb^2 \E_n^{\otimes 2} \left[ \gd_i\gd_j
            e^{H_{n,h}^{\go,(1)}+H_{n,h}^{\go,(2)}} e^{\mu \gk_n ((S_n^{(1)})^2+(S_n^{(2)})^2)}\right]
            \left(Z_{n,h}^{\mu ,(1)}Z_{n,h}^{\mu,(2)}\right)^{\gga-1}.
\end{multline}
Therefore,
\begin{multline}
 \sum_{i=1}^{2^{n}}\sum_{j=2^n+1}^{2^{n+1}}  \frac{\partial^2}{\partial \go_i \partial\go_j}
      \left(Z_{n,h}^{\mu ,(1)}Z_{n,h}^{\mu,(2)}\right)^{\gga} \\
 \leq \frac{\gamma^2\gb^2}2\E_n^{\otimes 2} \left[ \left( (S_n^{(1)})^2 + (S_n^{(2)})^2 \right)
            e^{H_{n,h}^{\go,(1)}+H_{n,h}^{\go,(2)}} e^{\mu \gk_n ((S_n^{(1)})^2+(S_n^{(2)})^2)}\right]
     \left(Z_{n,h}^{\mu ,(1)}Z_{n,h}^{\mu,(2)}\right)^{\gga-1},
\end{multline}
and as a consequence, since we chose $\gk_n=\gk^n$
\begin{equation}
 \frac{\partial \Phi}{ \partial t} (t,\mu) \leq \frac{\gk}2 \frac{\partial \Phi}{ \partial \mu} (t,\mu).
\end{equation} 
Thus, the function $t\mapsto \Phi (t, \mu-\gk t/2)$ is non-increasing and
\begin{equation}
  \log \Eo\left[ \left(Z_{n,h}^{\mu -\gk/2 ,(1)}Z_{n,h}^{\mu-\gk/2,(2)}\right)^{\gga} \right] =\Phi(1,\mu-\gk/2)
\leq \Phi(0,\mu) = 2\log \Eo_t\left[ \left(Z_{n,h}^{\mu } \right)^{\gga} \right].
\end{equation} 
Then, one uses that for $\mu \geq \frac{\gk}{2(1-2\gk)}$ one has $2\mu\gk \leq \mu-\gk/2$, which allows us to conclude.
\end{proof}

\subsection{Change of measure }
\label{sec:changmeas}

In this section we prove the lower bounds of Theorem \ref{thm:shift} on the
critical point shift for $B\le B_c$.


One fixes $\gga$ close to $1$ such that $x_{\gga}$ is also close to $1$, and
proves that if $h=\hca+u$ with $u>0$ small enough, one has $A_{n_0}:=\Eo\left[ \left( \bar Z^\go_{n_0,h}\right)\right] \leq x_{\gga}$ for some  
$n_0\in\N$.
To this purpose, we introduce a change of measure in the spirit of \cite{GLT09}.
Define
\begin{equation}
\label{def:gF}
 \begin{split}
   g(\go)& := \ind_{\{F(\go)\leq R\}} + \gep_R \ind_{\{F(\go)> R\}}, \\
 F(\go) &:= \langle V\go , \go \rangle - \Eo[\langle V\go , \go \rangle],
 \end{split}
\end{equation}
where the choices of the symmetric $2^n\times 2^n$ matrix $V$, of $R\in\R$ 
and $\gep_R>0$ will be made later. Note that we have chosen $F$ to be centered.
Then
using the H\"older inequality, one has
\begin{equation}
  \Eo\left[ (\bar Z^\go_{n,h})^{\gga}\right] = \Eo\left[g(\go)^{-\gga} (g(\go)\bar Z^\go_{n,h})^{\gga}\right]\leq
  \Eo\left[ (g(\go))^{-\frac{\gga}{1-\gga}}\right]^{1-\gga} \Eo\left[ g(\go) \bar Z^\go_{n,h}\right]^{\gga}.
\label{eq:holder}
\end{equation}

\begin{rem}\rm
  The original idea \cite{GLT08} is to take $g(\go) = \frac{\dd
    \check\Po}{\dd \Po}$ where $\check\Po$ is a new probability
  measure on $\{\go_1,\ldots,\go_{2^n}\}$ such that $\check\Po$ and
  $\Po$ are mutually absolutely continuous.  Then, to control both
  terms in \eqref{eq:holder}, one has to choose $\check\Po$ in a
  certain sense close enough to $\Po$, such that the first term is
  close to $1$, but also such that under the measure $\check\Po$ the
  annealed partition function $\Eo\left[ g(\go) \bar
    Z_{n,h_c^{\a}}\right]=\check\Eo\left[ \bar Z_{n,h_c^{\a}}\right]$
  is small.

  The choice of $g$ and $F$ in \eqref{def:gF} has the same effect of
  the change of measure in \cite{GLT08}, that is inducing negative
  correlations between different $\go_i$, and the specific form
  \eqref{def:gF} is chosen for technical reasons, to deal more easily
  with the case in which $\langle V\go , \go \rangle$ is large.
\end{rem}

Let us first deal with the Radon-Nikodym part of \eqref{eq:holder}: we make here the choice
$\gep_R:= \Po(F(\go)\geq R)^{1-\gga}$. Then one has
\begin{equation}
\label{eq:firstpart}
 \Eo\left[ (g(\go))^{-\frac{\gga}{1-\gga}}\right]  \leq  1 + (\gep_R)^{-\frac{\gga}{1-\gga}} \Po(F(\go)\geq R)
    =  1+ \Po(F(\go)\geq R)^{ 1-\gga} = 1+\gep_R.  
\end{equation}

We now use the following lemma to estimate $\gep_R$ in terms of $R$. We let $\|V\|^2=
\sum_{i,j}V_{ij}^2$ and $K$ denote the covariance matrix $(\kappa_{ij})_{1\le
i,j\le 2^n}$.
\begin{lemma}
\label{lem:normV}
If $V$ is such that $V_{ij}$ depends only on $d(i,j)$ and $\|V\|^2=1$,
then one has $\Var(F)<2 K_{\infty}^2 $ with $K_\infty$ defined in \eqref{eq:assumpK}, so that
\begin{equation}
 \Po(F(\go)\geq R) \leq \frac{2K_{\infty}}{R^2}\stackrel{R\to\infty}{\longrightarrow} 0 .
\end{equation}
Thus one gets that $\gep_R\leq const\times R^{-2(1-\gga)}$, which can be made arbitrarily small choosing $R$ large.
\end{lemma}

\begin{proof}
We have that $\Var(F) = \Eo\left[\langle V\go , \go \rangle^2 \right] - \Eo\left[\langle V\go , \go \rangle\right]^2$,
and we can compute
\begin{eqnarray}
 \Eo\left[\langle V\go , \go \rangle^2 \right] &=&
     \sum_{i,j=1}^{2^n} \sum_{k,l=1}^{2^n} V_{ij} V_{kl} \Eo[\go_i\go_j\go_k\go_l] 
    =\sum_{i,j=1}^{2^n} \sum_{k,l=1}^{2^n} V_{ij} V_{kl} (\gk_{ij} \gk_{kl} + \gk_{ik} \gk_{jl} + \gk_{il} \gk_{jk}) \nonumber \\
  & = & \Eo\left[\langle V\go , \go \rangle\right]^2 + 2 Tr\left( (VK)^2 \right).
\end{eqnarray}
We now use Lemma \ref{lem:matrix}, which says that $V$ and $K$ can be
codiagonalized, and that the eigenvalues of $K$ are bounded by
$K_{\infty}$, to get that $Tr\left( (VK)^2 \right) \leq K_{\infty}^2
Tr(V^2)=K_{\infty}^2$ (recall that $Tr(V^2)= \|V\|^2=1$, as $V$ is
symmetric). One finally gets that $\Var(F) \leq 2 K_{\infty}^2$, and
as $F$ is centered, using Tchebyshev's inequality gives the
result.
\end{proof}

Next,  we study the second factor in the r.h.s. of \eqref{eq:holder}:
\begin{equation}
 \Eo\left[ g(\go) \bar Z_{n,h}^{\go}\right] \leq
    \Eo\left[ \ind_{\{F(\go)\leq R\}} \bar Z_{n,h}^{\go}\right] +\gep_R \Eo\left[ \bar Z_{n,h}^{\go}\right].
\label{eq:plusqu1}
\end{equation}
To study the first term 
we define the measure $\Pot$ on $\{\go_1,\ldots,\go_{2^n}\}$ to be absolutely continuous
with respect to $\Po$, with Radon-Nikodym derivative given by
$
\frac{\dd \Pot}{\dd \Po} = \frac{\bar Z_{n,h}^{\go}}{\bar Z_{n,h}^{\a}}.
$
One then has
\begin{equation}
 \Eo\left[ \ind_{\{F(\go)\leq R\}} \bar Z_{n,h}^{\gb,\go}\right] = \bar Z_{n,h}^{\a} \Pot \left( F(\go)\leq R \right).
\label{eq:boundPot}
\end{equation}

\medskip
We are now ready to choose $V=V_n$, and we do so as in \cite{GLT08}.
We take $V$ to be zero on the diagonal ($V_{ii}=0$), and for $i,j \in\{1,\ldots,2^n\}$
\begin{equation}
\label{eq:V}
 V_{ij}:= \frac{\bE_n[\gd_i\gd_j]}{Y_n}, \indent \text{ if }  i\neq j , 
\end{equation}
where
\begin{equation}
\label{eq:Y}
 Y_n := \left(\sumtwo{i,j=1}{i\neq j}^{2^n} \bE_n[\gd_i\gd_j]^2 \right) ^{1/2}
\end{equation}
is used to normalize $V$. We stress that $V$ satisfy the conditions of Lemma \ref{lem:normV}.

 One can compute easily $Y_n$ , since from Proposition \ref{prop:nodi} we have $\bE_n[\gd_i\gd_j]= B^{-n-d(i,j)+1}$, and one finds (cf. \cite[Eq. (8.23)]{GLT08})
\begin{equation}\label{eq:calculY} 
Y_n=\begin{cases}
\sqrt{n} &\indent \text{ if } B=B_c:=\sqrt{2},\\
\Theta\left( \left(\frac{2}{B^2} \right)^n \right) &\indent \text{ if } B<B_c
\end{cases}
\end{equation}
where $X=\Theta(Y)$ means
that
$X\ge c Y$ for some positive constant $c$.

\begin{proposition}
\label{prop:PotF}
 We choose $V=V_n$ as in \eqref{eq:V}-\eqref{eq:Y}, and $R=R_n:= \frac12 \Eot[F(\go)]$.
Then there exists some $\gd>0$ small such that, if $u(2/B)^n\leq \gd$, one has
\begin{equation}
 R := \frac12 \Eot[F(\go)]\ge c  \gb^2 Y_n.
\label{eq:EotF}
\end{equation}
Therefore,  from \eqref{eq:calculY}, $R$ can be made arbitrarily large with $n$.
Moreover there exists a constant $\gz>0$ which does not depend on $n$, such that
\begin{equation}
 \Pot\left(F(\go) \geq R \right) =\Pot\left(F(\go) \geq \frac12 \Eot[F(\go)] \right)\geq \gz .
\label{eq:mainprop}
\end{equation}
\end{proposition}
Combining this Proposition to \eqref{eq:plusqu1} and \eqref{eq:boundPot}, one gets that
\begin{equation}
 \Eo \left[ g(\go) \bar Z_{n,h}^{\go}\right] \leq
        \bar Z_{n,h}^{\a} \left( 1-\gz + \gep_R \right).
\label{eq:secondpart}
\end{equation}

Recalling the 
equality \eqref{eq:thmbar} (which is the analog of Theorem \ref{thm:boundEnc} for
the alternative partition function $\bar Z_{n,h}^{\a}$),
one has for $\gk<B^2/4\wedge 1/2$
\begin{equation}
 \bar Z_{n,h}^{\a} \leq \E_n\left[ e^{c_1' u S_n} \right]+c_2'\gb^2\left(\frac{4\gk}{B^2}\right)^n\leq e^{c\gd}+\gd,
\end{equation} 
provided that $u\leq \gd (B/2)^n$ with $\gd$ small (to be able to apply Lemma \ref{lem:estimEsum} to $\E_n\left[ e^{c_1' u S_n} \right]$),
and that $n\geq n_{\gd}$ to deal with the term $\left(4\gk/B^2\right)^n$.
Therefore, if $\gd$ and $\gep_R$ was chosen small enough (that is smaller than some constant $c=c(\gz)$),
one has for $n\geq n_{\gd}$ that $\Eo \left[ g(\go) \bar Z_{n,h}^{\gb,\go}\right]\leq 1-\gz/2$ for all $u\leq \gd (B/2)^n$.
This and \eqref{eq:firstpart} bound the two terms in \eqref{eq:holder}, so that one has
\begin{equation}
A_n:=\Eo\left[ (\bar Z_{n,h}^{\go})^{\gga}\right] \leq \left( 1+\gep_R \right)
   \left( 1-\gz/2\right)^{\gga}\leq 1-\gz/3\leq x_{\gga},
\end{equation}
where the two last inequalities hold if $\gep_R$ is small and $\gga$ close to $1$.
To sum up, for $\gd,\gb$ small and $R$ large enough,
one has that $A_n\leq x_{\gga}$
for all $u\leq \gd (B/2)^n$, and so $\F(\gb,\hca+u)=0$.

Then, let us check how large has to be $n$ so that our choice of $R:=\frac12 \Eot[F(\go)]$
becomes large. From 
Proposition \ref{prop:PotF} one has that $R \geq c \gb^2 Y_n$ so that one has to take
$\gb^2 Y_n\geq C$ for some constant $C$ large enough.
From  \eqref{eq:calculY}, in order to have $\gb^2 Y_n\geq C$,
\begin{itemize}
 \item if $B<B_c$, it is enough to take
   $n$ larger than $n_0:= \log(C'\gb^{-2})/\log(2/B^2)$;

 \item if $B=B_c$, one has to take
   $n$ larger than $n_0:= c'\gb^{-4}$. 
\end{itemize}

Then for $n= n_0$ one gets that $R$ is large,
but one also needs to take $u\leq \gd(2/B)^{n_0}$ to ensure that $A_{n_0\vee n_{\gd}}\leq x_{\gga}$ .
Notice that from the choice of $n_0$ above, the condition on $u$ translates into
\begin{equation}
u \leq 
 \begin{cases}
  c' \gb^{2 \log(2/B)/ \log(2/B^2)} =c' \gb^{\frac{2}{2-\nu}}& \indent \text{if } B<B_c,\\
  e^{-c \gb^{-4}} & \indent \text{if } B=B_c,
 \end{cases}
\end{equation}
where we also used that $\nu=\log 2/ \log(2/B)$. One then gets the
desired bounds \eqref{alb}-\eqref{ulb} on the difference between
quenched and annealed critical points.
\qed

\subsection{Proof of Proposition \ref{prop:PotF}}
To compute $\Eot[F(\go)]$,
we define for any $1\leq i,j\leq 2^n$
\begin{equation}
 U_{ij}:=\Eot[\go_i \go_j] = \frac{1}{\bar Z_{n,h}^{\a}}\bE_n \Eo \left[ \go_i \go_j \, e^{ \bar H_{n,h}^{\go}}\right].
\end{equation}
A Gaussian integration by parts gives easily
\begin{eqnarray}
\label{eq:def u_ij}
 U_{ij}  
=   \gk_{ik}+u_{ij}:=  \gk_{ij} + \gb^2 \sum_{k,l=1}^{2^n} \gk_{ik}\gk_{jl} \bar \bE_{n,h}^{\a} [\gd_k \gd_l],
\end{eqnarray}
where $\bar \bE_{n,h}^{\a}$ denotes expectation w.r.t. the measure
whose density with respect to $\bP_n$ is $\exp(\bar H^{\a}_{n,h})/\bar
Z^{\a}_{n,h}$.  We then compare $\bar \bE_{n,h}^{\a} [\gd_k \gd_l]$
with $\bE_n [\gd_k \gd_l]$, using that $h=\hca+u$, $0\leq u\leq
\gd(B/2)^n$:
\begin{equation}
 \bar \bE_{n,h}^{\a} [\gd_k \gd_l] = \frac{1}{\bar Z_{n,h}^{\a}}
        \bE_n\left[ \gd_k \gd_l e^{\bar H_{n,\hca}^{\a}} e^{u S_n}\right]
 \leq e^{2c_1\gb^2} \bE_n\left[ \gd_k \gd_l e^{ e^{c_1\gb^2} u S_n}\right]
  \leq c' \bE_n\left[ \gd_k \gd_l \right]
\label{useProp+Cor}
\end{equation} 
where in the first inequality we used Remark \ref{remarca} and also
the fact that $\bar Z_{n,h}^{\a}\geq \bar Z_{n,\hca}^{\a}\geq 1$, and
in the second inequality we used that $u(2/B)^n\leq \gd$ to apply
Corollary~\ref{cor:estimsum}.  The same argument easily gives $\bar
\bE_{n,h}^{\a} [\gd_k \gd_l]\geq c \bE_n\left[ \gd_k \gd_l \right]$ in
the range of $u$ considered, so that $ c\gb^2 a_{ij}\leq u_{ij}\leq
c'\gb^2 a_{ij}$, where
\begin{equation}
\label{eq:def a_ij}
a_{ij}:=\sum_{k,l=1}^{2^n} \gk_{ik}\gk_{jl} \bE_{n} [\gd_k \gd_l]  \ge Y_n (KVK)_{ij}
\end{equation}
(the inequality is due to the fact that $V$ is zero on the diagonal).
We finally get
\begin{equation}
 \Eot[F(\go)]=\Eot[\langle V\go, \go \rangle] - \Eo[\langle V \go,\go \rangle]
    = \sum_{i,j=1}^{2^n} V_{ij} (\gk_{ij} +  u_{ij}) - \Eo[\langle V \go,\go \rangle]
      = \sum_{i,j=1}^{2^n} V_{ij}u_{ij},
\end{equation}
so that we only have to compute
$\sum_{i,j=1}^{2^n} V_{ij}a_{ij}\ge Y_n {\rm Tr}(VKVK)
.
$ 
Since $\|V\|^2=1$ and all eigenvalues of $K$
are between $1$ and $K_{\infty}$, one has
${\rm Tr}\left( (VK)^2 \right) = \Theta(1)$.
Altogether, we get \eqref{eq:EotF}.

\medskip

We now prove \eqref{eq:mainprop}. Using the Paley-Zygmund inequality,
we get that
\begin{equation}
\Pot(F(\go)\geq R) = \Pot \left( F(\go)\geq \frac12 \Eot[F(\go)] \right) \geq \frac{\Eot[F(\go)]^2}{4 \Eot[F(\go)^2]},
\end{equation}
so that we only have to prove the following:
\begin{equation}
 \tilde{\Var}(F(\go)) = \Eot[\langle V\go,\go \rangle^2 ]- \Eot[\langle V\go,\go \rangle ]^2=O(\Eot[F(\go)]^2).
\label{eq:varF}
\end{equation}
Indeed from this it follows immediately that there exists some constant $\gz>0$
such that $\Eot[F(\go)]^2 /\Eot[F(\go)^2]\geq \gz$ . 

We now prove \eqref{eq:varF}, studying
$\Eot\left[\langle V\go,\go \rangle^2\right] = \sum_{i,j,k,l=1}^{2^n} V_{ij} V_{kl}  \Eot [\go_i\go_j\go_k\go_l]$,
starting with the computation,
for any $1\leq i,j,k,l\leq 2^n$, of
\begin{equation}
 \Eot [\go_i\go_j\go_k\go_l] =\frac{1}{\bar Z_{n,h}^{\a}}\bE_n \Eo \left[ \go_i \go_j \go_k\go_l\ 
   e^{ \bar H_{n,h}^{\go}}\right].
\end{equation}
Again, a Gaussian integration by parts gives, after elementary computations,
\begin{multline}
\Eot \left[ \go_i \go_j \go_k \go_l\right] =A_{ijkl}+B_{ijkl}
:=  \left[\gk_{ij} U_{kl}    +\gk_{ik}  U_{jl}    +\gk_{il}  U_{jk} + \gk_{jk}u_{il} +\gk_{jl} u_{ik} + \gk_{kl}u_{ij}\right]\\
   + \gb^4 \sum_{r,s,t,v=1}^{2^n} \gk_{ir}\gk_{js}\gk_{kt}\gk_{lv} 
     \bar \bE_{n,h}^{\a}\left[ \gd_r\gd_s\gd_t\gd_v  \right].
\label{eq:omegaijkl}
\end{multline}
We estimate $\Eot\left[\langle V\go,\go \rangle^2\right]$ by analyzing
separately $A_{ijkl}$ and $B_{ijkl}$.

\medskip

\emph{Contribution from $B_{ijkl}$:} we have
\begin{eqnarray}
 B_{ijkl} 
  \leq  c \beta^4\sum_{r,s,t,v=1}^{2^n} \gk_{ir}\gk_{js}\gk_{kt}\gk_{lv} \bE_n\left[ \gd_r\gd_s\gd_t\gd_v  \right],
\end{eqnarray}
where we used again Proposition \ref{prop:Pnc} and Corollary \ref{cor:estimsum} as in \eqref{useProp+Cor} (recall
that we consider $u\leq \gd (B/2)^n$). Then defining
\begin{eqnarray}
  \label{eq:8}
  W_{ij}:=\frac{\bE_n[\delta_i\delta_j]}{Y_n}=V_{ij}+\frac{{\bf 1}_{\{i=j\}}}{Y_n B^n}.
\end{eqnarray}
we get
\begin{multline}
 \sum_{i,j,k,l=1}^{2^n} V_{ij}V_{kl} B_{ijkl} \leq
          c \beta^4\sum_{r,s,t,v=1}^{2^n} (KWK)_{rs} (KWK)_{tv} \bE_n\left[ \gd_r\gd_s\gd_t\gd_v  \right] \\
   \leq c'\beta^4 \sumtwo{r,s,t,v=1}{r\neq s, t\neq v}^{2^n} W_{rs} W_{tv} \bE_n\left[ \gd_r\gd_s\gd_t\gd_v  \right]
    + c''\beta^4 \sum_{r,t,v=1}^{2^n}W_{rr}W_{tv}\bE_n\left[ \gd_r\gd_t\gd_v  \right],
\label{eq:crossedterm}
\end{multline}
where we used the following claim:
\begin{claim}
There exists a constant $c'>0$ such that for every $1\leq i,j\leq 2^n$, $(WK)_{ij}\leq c'W_{ij}$
and $(KW)_{ij}\leq c'W_{ij}$.
\label{claim:KV}
\end{claim}

\begin{proof}[Proof of the Claim]
We write $q=d(i,j)$, so $W_{ij}=:W_q$, and
\begin{equation}
\label{eq:VK}
(WK)_{ij} = \sum_{l=1}^{2^n} W_{il}\gk_{lj} = \sum_{p=0}^{q-1} 2^{p-1} W_p \gk_q + 
       \sum_{p=0}^{q-1} 2^{p-1} W_q\gk_p + \sum_{p=q+1}^{n} 2^{p-1} W_p \gk_p,
\end{equation}
where we decomposed the sum according to the positions of $l$ ($d(i,l)=p<q$, $d(i,l)=q$ or $d(i,l)>q$).
Using that $W_p$ is decreasing with $p$,
we get that the second and the third term are
both smaller than $(\sum 2^p \gk_p) W_q$. We only have to deal with the first term, using the explicit expression of
$W_p$,
together with Proposition \ref{prop:nodi}:
\begin{equation}
 \sum_{p=0}^{q-1} 2^{p-1} W_p = \frac{1}{Y_n} B^{-n}  \sum_{p=0}^{q-1} \left( \frac{2}{B} \right)^{p-1}
    \leq c\, \frac{1}{Y_n} B^{-n} \left( \frac{2}{B} \right)^{q} = c\, 2^q W_q,
\end{equation}
so that the first term in \eqref{eq:VK} is smaller than $c 2^q \gk_q W_q$. One then has that $(WK)_{ij}\leq c' W_{ij}$,
and the same computations also gives that $(KW)_{ij}\leq c' W_{ij}$.
\end{proof}

The main term in the r.h.s. of \eqref{eq:crossedterm} is the first
one, for which we have 
\begin{lemma}
\label{lem:crossedterms}
 Let $B\leq B_c$. There exists a constant $c>0$ such that
  \begin{multline}
  \sumtwo{r,s,t,v=1}{r\neq s, t\neq v}^{2^n} V_{rs} V_{tv} \bE_n\left[
    \gd_r\gd_s\gd_t\gd_v  \right]= \frac{1}{Y_n^2}\sumtwo{r,s,t,v=1}{r\neq s, t\neq v}^{2^n} 
            \E_n[\gd_r\gd_s] \E_n[\gd_t\gd_v] \E_n[\gd_r\gd_s\gd_t\gd_v]
       \leq c Y_n^2.
  \end{multline}
\end{lemma}
This can be found in the proof of Lemma 4.4 of \cite{GLT08} for $B=B_c$;
the proof is easily extended to the case $B<B_c$.

As for the remaining terms in \eqref{eq:crossedterm}, it is not hard
to see, using repeatedly Proposition~\ref{prop:nodi},  that they give
a contribution of order $o(Y_n^2)$.
For instance, one has
\begin{eqnarray}
  \beta^4 \sumtwo{r,t, v=1}{t\ne v}^{2^n}W_{rr}W_{tv}\bE_n\left[ \gd_r\gd_t\gd_v  \right]\le c\beta^4\frac{2^{n}}{B^n Y_n^2}\sum_{p=0}^n
2^p B^{-n-p}\sum_{q=0}^n2^q B^{-n-p-q}=\gb^4 o(Y_n^2).
\end{eqnarray}
Altogether one has
\begin{equation}
\sum_{i,j,k,l=1}^{2^n} V_{ij}V_{kl} B_{ijkl} = \gb^4 O\left(Y_n^2\right)=O\left(\Eot[F(\go)]^2 \right),
\label{alto}
\end{equation}
cf. \eqref{eq:EotF}.
\medskip

\emph{Contribution of $A_{ijkl}$:}
recalling that $U_{ij}=\gk_{ij} + u_{ij}$, we have $\gk_{ij} U_{kl}+ \gk_{kl}u_{ij} \leq U_{ij}U_{kl}$. Thus, we get
\begin{equation}
\label{basso} \sum_{i,j,k,l=1}^{2^n} V_{ij}V_{kl} (\gk_{ij} U_{kl}+ \gk_{kl}u_{ij}) \leq \left( \sum_{i,j=1}^{2^n} V_{ij}U_{ij} \right)^2
   = \Eot\left[\langle V\go,\go\rangle \right]^2,
\end{equation}
that we recall is not $O(\Eot[F(\go)]^2)$, but will be canceled in the variance.
The other contributions are, thanks to symmetry of $V$, all equal to (or smaller than)
\begin{equation}
 \sum_{i,j,k,l=1}^{2^n} V_{ij}V_{kl} \gk_{ik}U_{jl} = \sum_{i,j,k,l=1}^{2^n} V_{ij}V_{kl} \gk_{ik}\gk_{jl} +
   \sum_{i,j,k,l=1}^{2^n} V_{ij}V_{kl} \gk_{ik}u_{jl},
\end{equation}
where the first term is ${\rm Tr}\left((VK)^2\right)$ which is bounded as remarked before.
Thanks to the estimate $u_{jl} \leq c'\gb^2 a_{jl}= c'\gb^2 Y_n(KWK)_{jl} $, the second term is bounded
above by a constant times 
\begin{gather}
\nonumber
\beta^2
Y_n \sum_{i,j,k,l=1}^{2^n} V_{ij}V_{kl} \gk_{ik} (KWK)_{jl} \le \beta^2 Y_n {\rm Tr} \left( (WK)^3 \right) \\
\label{medio}
\leq  c\beta^2 Y_n {\rm Tr} (W^2) \le 2 c\beta^2 Y_n=O(\Eot[F(\go)]).
\end{gather}
We used Lemma \ref{lem:matrix} to codiagonalize $W$ and $K$
and to bound the eigenvalues of $K$ by a constant, and then the fact that
the eigenvalues $\lambda_i$ of $W$ are also bounded, so that
$\sum |\lambda_i|^3\le c\sum |\lambda_i|^2=c{\rm Tr}(W^2)=O(1)$.
Indeed, ${\rm Tr}(W^2)={\rm Tr}(V^2)+\sum_i W^2_{ii}=1+(2/B^2)^nY_n^{-2}=1+o(1)$.
Putting together \eqref{eq:omegaijkl} with the estimates \eqref{alto},
\eqref{basso} and \eqref{medio}  we have
\begin{gather}\nonumber
  \tilde{\Var}(F(\go)) = \tilde{\mathbb E} (\langle
  V\go,\go\rangle^2)-\left( \tilde{\mathbb E} \langle
    V\go,\go\rangle\right)^2=
  \sum_{ijkl}(A_{ijkl}+B_{ijkl})V_{ij}V_{kl}-\left( \tilde{\mathbb E}
    \langle V\go,\go\rangle\right)^2\\= O\left(\Eot[F(\go)]^2 \right)
\end{gather}
and \eqref{eq:varF} is proven.

\begin{appendix}

\section{Pure model estimates}
\label{sec:purestim}

We first give some estimates on the partition function of a system of size $n$.
\begin{lemma}
\begin{enumerate}
\item There exist constants $a_0>0$ and $c_0>0$ such that for any $n\geq 0$, if $u\leq a_0\, (B/2)^n$ one has
\begin{equation}
\E_n\left[\exp \left(u S_n\right)\right] \leq \exp\left(c_0 u (2/B)^n\right).
\end{equation}
\item There exists  a constant $c>0$ such that for any $n\geq 0$ and $u\geq 0$ one has
\begin{equation}
 \E_n \left[ \exp \left( u S_n \right) \right]  \leq c\exp\left( c u^{\nu}2^n \right),
\end{equation}
where $\nu$ is as in \eqref{eq:purenu}.
\end{enumerate}
\label{lem:estimEsum}
\end{lemma}

\begin{proof}
For the first inequality, the same type of computation was already done in \cite{GLT07}, and we give here only an outline of the proof.
The partition function $R_k$ of the pure model satisfies  the iteration 
\begin{equation}
 \begin{cases}
  R_0 = e^u, \\
  R_{k+1} = \frac{R_k^2+B-1}{B}.
 \end{cases}
\end{equation}
Defining $P_k:=R_k-1$
it is easy to show by recurrence that $P_k \leq c_0 u \left( \frac{2}{B} \right)^k$ for every $k\leq n$
(because we stay in the linear regime for the chosen value of $u$),
so that for $k=n$ we get the result. 

For the second inequality, we use that for any $n\geq 0$ and $u\geq 0$,
\begin{equation}
 \frac{1}{2^n} \log \E_n \left[ \exp \left( uS_n\right) \right] \leq \F(u)+ \frac{c(B)}{2^n},
\end{equation}
from \cite[Theorem 1.1]{GLT07}, and this gives immediately the result, using 
\eqref{eq:pureF}.
\end{proof}

Defining for any subset $I\subset \{1,\ldots,2^n\}$
$\gd_I:=\prod_{i\in I} \gd_i$, and $\gd_I=1$ if $I=\emptyset$, one
wants to compare $\E_n[\gd_I e^{uS_n}]$ and $\E_n[\gd_I]$ when the
partition function $Z_{n,h}^{\rm pure}$ is still in the linear regime $0\le u\leq
a_0\ (B/2)^n$, the bound $\E_n[\gd_I e^{uS_n}]\geq \E_n[\gd_I]$ being
trivial.
\begin{cor}
\label{cor:estimsum}
There exist constants $a_0>0$ and $c'>0$ such that for any $n\geq 0$
and any non-empty subset $I\subset \{1,\ldots,2^n\}$, if $0\le u\leq a_0\, (B/2)^n$ one has
\begin{equation}
\E_n\left[\gd_I \exp \left(u S_n\right)\right] \leq
  \exp\left(c' u \left(\frac{2}{B}\right)^n\right)^{|I|} \E_n\left[ \gd_I\right].
\end{equation}
\end{cor}

\begin{proof}
 We prove  by iteration on $n$ that for all non-empty subsets $I\subset \{1,\ldots,2^n\}$, if $u\leq a_0\, (B/2)^n$ one has
\begin{equation}
\E_n\left[\gd_I \exp \left(u S_n\right)\right] \leq
  \exp\left(c_0 u \sum_{k=0}^{n} \left(\frac{2}{B}\right)^k\right)^{|I|} \E_n\left[ \gd_I\right],
\end{equation}
where $c_0$ is the constant obtained in Lemma \ref{lem:estimEsum}.

The case $n=0$ is trivial.
Let us assume that we have the assumption for $n\geq0$ and prove it for $n+1$.
Take
$I$ a non-empty subset of $\{1,\ldots,2^{n+1}\}$.
As in the proof of Lemma
\ref{presquepropPnc}, one decomposes $I$ into its ``left'' and
``right'' part and
writes $\E_{n+1}[\gd_{I}]=
\frac{1}{B}\E_{n}[\gd_{I_1}]\E_{n}[\gd_{\tilde I_2}]$ and
$|I|=|I_1|+|\tilde I_2|$.

If $I_1,\tilde I_2\neq\emptyset$, using the induction hypothesis, 
one easily has
\begin{multline}
 \E_{n+1}\left[\gd_I \exp \left(u S_{n+1}\right)\right]
  = \frac{1}{B} \E_{n}\left[\gd_{I_1} \exp \left(u
      S_{n}\right)\right]\E_{n}\left[\gd_{\tilde I_2} \exp \left(u S_{n}\right)\right] \\
  \leq \exp\left(c_0 u \sum_{k=0}^{n}
    \left(\frac{2}{B}\right)^k\right)^{|I_1|+|\tilde I_2|}
  \frac{1}{B}\E_{n}[\gd_{I_1}]\E_{n}[\gd_{\tilde I_2}],
\end{multline} 
which gives the right bound.

If $I_1= \emptyset$ (or analogously if $\tilde I_2=\emptyset$), one has
$\E_{n+1}[\gd_{I}]= \frac{1}{B}\E_{n}[\gd_{\tilde I_2}]$ and
\begin{multline}
 \E_{n+1}\left[\gd_I \exp \left(u S_{n+1}\right)\right]
  = \frac{1}{B} \E_{n}\left[\exp \left(u
      S_{n}\right)\right]\E_{n}\left[\gd_{\tilde I_2} \exp \left(u S_{n}\right)\right] \\
  \leq e^{c_0 u(2/B)^{n+1}} \exp\left(c_0 u \sum_{k=0}^{n}
    \left(\frac{2}{B}\right)^k\right)^{|\tilde I_2|}
  \frac{1}{B}\E_{n}[\gd_{\tilde I_2}],
\end{multline} 
where the first part is dealt with Lemma \ref{lem:estimEsum}, and the second one with the induction hypothesis.
\end{proof}

\begin{theorem}
\label{thm:bound2}
Let $B\in(1,2)$.
 Let $(b_{n})_{n\geq 0}$ be a sequence that goes to $0$ as $n$ goes to infinity. There exists a constant $c_b>0$ such
that for all $n\geq 0$ and every $0\leq u\leq b_n(\frac{B^2}{4} \wedge \frac12 )^n$ one has
\begin{equation}
\E_n\left[  \exp\left( u (S_n)^2 \right)  \right]
     \leq \exp\left( c_b u\left(\frac{4}{B^2}  \right)^n \right).
\end{equation}
\end{theorem}

\begin{cor}
\label{cor:bound2}
Let $B\in(1,2)$,
$\gk<\frac{B^2}{4} \wedge \frac12$ and note $\gp:=(2\gk)\vee \frac{4\gk}{B^2}<1$. Then for every $A>0$ there exists a constant $c_A>0$
such that for any $n\ge 0$, any $u\in[0,A]$ and any subset $I$ of $\{1,\ldots,2^n\}$, one has
\begin{equation}
\E_n \left[\gd_I \exp\left(u\gk^n (S_n)^2 \right) \right] \leq \left( e^{c_A u\gp^n }\right)^{n|I|+1} \E_n[\gd_I]. 
\end{equation}
Note that if $I=\emptyset$, the statement is implied by Theorem \ref{thm:bound2}.
\end{cor}



\begin{proof}[Proof of Theorem \ref{thm:bound2}]

The proof relies on Lemma \ref{lem:estimEsum}.
Consider $u\leq b_n(\frac{B^2}{4} \wedge \frac12)^n$.
One writes 
\begin{equation}
J:= \E_n\left[\exp\left(\frac12 u (S_n)^2\right)\right]=\frac{1}{\sqrt{2\pi}}\int_{-\infty}^{+\infty} e^{-z^2/2} \E_n\left[\exp\left(z\sqrt{u} S_n\right)\right]\dd z.
\end{equation}
 
One sets $\gD:=  \frac{a}{\sqrt{u}} \left(\frac{B}{2}\right)^n$, where $a$
is a constant that will be chosen small.
Note that thanks to our choice of $u$, one has
$\gD \geq a\, b_n^{-1/2}$ that goes to infinity as $n$ grows to infinity. Then one decomposes the
integral $J$ according to the values of $z$, and writes $J=J_1+J_2$, where
\begin{equation}
\begin{split}
J_1:= \frac{1}{\sqrt{2\pi}}\int_{z\leq \gD} e^{-z^2/2} \E_n\left[\exp\left(z\sqrt{u} S_n\right)\right] \dd z\\
J_2:=\frac{1}{\sqrt{2\pi}}\int_{z\geq \gD}e^{-z^2/2} \E_n\left[\exp\left(z\sqrt{u} S_n\right)\right] \dd z.
\end{split}
\end{equation}
To bound $J_1$, one chooses $a\leq a_0$ with $a_0$ as in Lemma \ref{lem:estimEsum}, such that for the values of $z$ considered one has $z\sqrt{u}\leq a_0 (B/2)^n$
and then one applies Lemma \ref{lem:estimEsum}-(1) to get
\begin{equation}
J_1 \leq \frac{1}{\sqrt{2\pi}}\int_{z\leq \gD} e^{-z^2/2} \exp\left(c z\sqrt{u} (2/B)^n\right)\dd z 
     \leq \exp \left( \frac{c^2}{2} u \left(4/B^2\right)^n\right).
\end{equation}

We deal with the term $J_2$, decomposing again according to the values of $z$. Let us first introduce some notations:
we define the sequence $(\gD_k)_{k\geq 0}$ by the iteration
\begin{equation}
 \begin{cases}
  \gD_0 = \gD \\
  \gD_{k+1}= \gD(\gD_k)^{2/\nu} \, (>\gD_k>1),
 \end{cases}
\end{equation}
and define also $m=\inf \{ k,\, \gD_k \geq A\sqrt{u}2^n \}$, for
some $A$ chosen large enough later. We point out that $m$ is finite.
Indeed for a fixed large $n$,
if $\nu\leq 2$, then $\gD_k \geq \gD^{k+1}$ and goes to infinity as $k$ goes to infinity.
Otherwise, if $\nu>2$, $\gD_k$ goes to $\gD^{\nu/(\nu-2)}$ as $k$ goes to infinity.
Then, we just need to check that $\gD^{\nu/(\nu-2)}\geq A\sqrt{u}2^n$ if $n$ is large. Using the value of
$\nu=\log 2/ \log(2/B)$ one has $2^{1/\nu}=2/B$, so that 
$\gD^{\nu}=a^{\nu} u^{-\nu/2}2^{-n}$. Then
\begin{equation}
 \frac{\gD^{\nu}}{(\sqrt{u}2^n)^{\nu-2}} = a^{\nu} \frac{u^{-\nu/2}2^{-n}}{u^{\nu/2-1}2^{n(\nu-2)}}
  =  a^{\nu} \left(u2^n\right)^{-(\nu-1)} \geq a^{\nu} b_n^{1-\nu},
\label{calculgDnu}
\end{equation}
where we used that $u2^n \leq b_n$. As $\nu>1$, it remains only to take $n$ large.

One decomposes $J_2$ as follows:
\begin{multline}
\label{eq:decompJ2}
 J_2 = \sum_{k=0}^{m-1} \frac{1}{\sqrt{2\pi}}\int_{\gD_k}^{\gD_{k+1}} e^{-z^2/2}
   \E_n \left[ \exp \left( z\sqrt{u}S_n \right) \right] \dd z\\
    + \frac{1}{\sqrt{2\pi}}\int_{ \gD_{m}}^{+\infty} e^{-z^2/2} \E_n \left[ \exp \left(  z\sqrt{u}S_n\right) \right] \dd z .
\end{multline}

Each term of the sum in \eqref{eq:decompJ2} can be dealt with Lemma \ref{lem:estimEsum}-(2). One gets
\begin{multline}
 \frac{1}{\sqrt{2\pi}}\int_{ \gD_k}^{\gD_{k+1}} e^{-z^2/2} \E_n \left[ \exp \left( z\sqrt{u}  S_n \right) \right] \dd z 
  \leq  \E_n \left[ \exp \left( \gD_{k+1}\sqrt{u}S_n \right) \right] P\left( \cN \geq \gD_{k} \right) \\
  \leq  c_1 \exp\left(c_ 2 2^n u^{\nu/2}(\gD_{k+1})^{\nu}\right) \exp\left( -c (\gD_k)^2\right),
\label{1boundJ2}
\end{multline}
where $\cN$ stands for a standard centered Gaussian.
Now recall the definition of $\gD_k$ and $\gD$, that gives
$(\gD_{k+1})^{\nu} = \gD^{\nu}(\gD_k)^2= a^{\nu} u^{-\nu/2}2^{-n}(\gD_k)^2$, so that one can bound the term in \eqref{1boundJ2} by
\begin{equation}
c_1 \exp\left((c_ 2a^{\nu}  -c) (\gD_k)^2\right)\leq c_1 \exp\left(-c (\gD_k)^2/2\right),
\end{equation}
where the inequality is valid provided one has chosen $a$ sufficiently small.

Let us now deal with the last term in \eqref{eq:decompJ2}, trivially
bounding $S_n\le 2^n$:
\begin{multline}
 \frac{1}{\sqrt{2\pi}}\int_{ \gD_{m}}^{\infty} e^{-z^2/2} \E_n \left[ \exp \left( z\sqrt{u}S_n \right) \right] \dd z
  \leq \frac{1}{\sqrt{2\pi}}\int_{ \gD_{m}}^{\infty} e^{-z^2/2} e^{z \sqrt{u}2^n} \dd z \\
  =  e^{u4^n/2} P\left( \cN \geq \gD_{m}-\sqrt{u}2^n \right) \leq e^{A^{-2}(\gD_m)^2}e^{-c (1-A^{-1})^2(\gD_{m})^2}\leq e^{-c(\gD_m)^2/2} ,
\end{multline}
where we used that $\sqrt{u}2^n \leq A^{-1}\gD_m$, and supposed that $A$ was chosen large enough for the last inequality.


We finally get that for $n$ large one has
\begin{equation}
 J_2 \leq c_1\sum_{k=0}^{m} e^{-c (\gD_k)^2/2} \leq 
\begin{cases}
 C e^{-c \gD^2/2}  &\text{ if } \nu\leq 2,\\
 C m e^{-c \gD^2/2}& \text{ if } \nu> 2,
\end{cases}
\end{equation}
where in the case $\nu\leq 2$ we used that $ \gD_k \geq \gD^{k+1}$.
Note that for $\nu>2$, using \eqref{calculgDnu}, one
also can  bound $m$ from above as follows: since
$\gD_k = \gD^{\frac{1-(2/\nu)^{k+1}}{1-2/\nu}}$,
\begin{equation}
 \frac{\gD_k}{\sqrt{u}2^n} =  \frac{\gD^{\nu/(\nu-2)}}{\sqrt{u}2^n} \gD^{-\frac{\nu}{\nu-2}(2/\nu)^{k+1} }
\geq a^{\nu/(\nu-2)} b_n^{(1-\nu)/(\nu-2)} \gD^{-c'(2/\nu)^{k} }.
\end{equation}
So if one takes $k\geq -\log \log \gD / \log(2/\nu)$ one gets that $\gD_k \geq a^{\nu/(\nu-2)} b_n^{(1-\nu)/(\nu-2)} e^{-c'} \sqrt{u}2^n$. 
If $n$ is large enough this implies that $m\leq const\times\log \log \gD$.

Then one easily gets that $J_2 = o\left(\gD^{-2}\right)$, with $\gD^{-2}=O\left(u(4/B^2)^n\right)$,
so that combining with the bound on $J_1$ one has
\begin{equation}
J\leq  \exp \left( \frac{c_0^2}{2}u \left(4/B^2\right)^n\right) + o\left( u \left(4/B^2\right)^n\right).
\end{equation}
\end{proof}

\begin{proof}[Proof of Corollary \ref{cor:bound2}]
We proceed by induction. Fix $A>0$ and $u\leq A$, and take the constant $c_A$ obtained in Theorem \ref{thm:bound2}
for the sequence $b_n = A \left( \frac{4\gk}{B^2}\wedge 2\gk\right)^{n}$.
The case $n=0$ is trivial. Suppose now that the assumption is true for
some $n$, and take $I$ a subset of $\{1,\ldots,2^{n+1}\}$.\\
Suppose $I\neq \emptyset$ (otherwise one already has the result from Theorem \ref{thm:bound2}).
As in the proof of Lemma
\ref{presquepropPnc}, one decomposes $I$ into its ``left'' and
``right'' part and $\E_{n+1}[\gd_{I}]=\frac{1}{B}
\E_{n}[\gd_{I_1}]\E_{n}[\gd_{\tilde I_2}]$.
Using that $(S_{n+1})^2\leq 2(S_n^{(1)})^2+2(S_n^{(2)})^2$ one gets
\begin{multline}
  \E_{n+1}\left[ \gd_I \exp\left( u \gk^{n+1} (S_{n+1})^2 \right)\right] \\
  \leq \frac{1}{B}\E_{n}\left[ \gd_{I_1} \exp\left( (2\gk) u \gk^n
      (S_{n})^2 \right)\right] \E_{n}\left[ \gd_{\tilde I_2} \exp\left( (2\gk) u \gk^n (S_{n})^2 \right)\right]\\
  \leq \frac{1}{B} \E_n[\gd_{I_1}] \E_n[\gd_{\tilde I_2}] \left( e^{c_A u (2\gk)
      \varphi^n} \right)^{n|I_1|+n|\tilde I_2|+2}
        \leq \E_{n+1} [\gd_I] \left( e^{c_A u 2\gk \varphi^n} \right)^{(n+1)|I|+1},
\end{multline}
where for the second inequality we used the recursion assumption and
for the last one the assumption $|I|\geq 1$.
Now one just uses that $2\gk\leq \gp$ to conclude.
\end{proof}


From Corollary \ref{cor:bound2} one can deduce the following Proposition,
useful to control the variance of the partition function (see Section \ref{sec:var}).
Define as in \eqref{defDn} $D_n:=\sum_{i,j=1}^{2^n} \gk_{ij} \gd_i
\gd'_j$,
where $\delta$ and $\delta'$ are the populations at generation $n$
of two independent GW trees.

\begin{proposition}
Let $B\in(1,2)$, $\gk< \frac12 \wedge \frac{B^2}{4}$ and set
$\varphi=(2\gk)\wedge(4\gk/B^2)<1$.

$\bullet$ If $B>B_c$, then for every $\Phi\in\left( \frac{2}{B^2}\vee \gp,1 \right)$
 there exist some $u_0>0$ and some constant $c>0$, such that for every $n\in\N$, $u\in [0,u_0]$ one has
\begin{equation}
  \E_n^{\otimes2} \left[\exp\left( uD_n \right)\right]
 \leq1+ cu \Phi^n.
\end{equation}

$\bullet$ If $B< B_c$
 there exist some $a_1>0$ and some constant $c>0$, such that for every $n\in\N$,
if $u\leq a_1 \left( \frac{B^2}{2} \right)^n$ one has
\begin{equation}
  \E_n^{\otimes2} \left[\exp\left( uD_n\right)\right]
 \leq 1+ cu \left(\frac{2}{B^2} \right)^n .
\end{equation}

$\bullet$ If $B=B_c$,
there exists some $u_0$ such that if $u\leq u_0$ then for all
$n \leq \frac12 u^{-1/3}$ one has
\begin{equation}
 \E_n^{\otimes2} \left[\exp\left( uD_n\right)\right]
 \leq 1+ 2 u^{1/3}.
\end{equation}
\label{prop:bound22}
\end{proposition}


\begin{proof}
One has
\begin{gather}
 D_{n+1}= D_n^{(1)}+D_{n}^{(2)} + \gk_{n+1} \left( S_{n}^{(1)}
   {S'}_{n}^{(2)}+ {S'}_{n}^{(1)}  {S'}_{n}^{(2)}\right) \nonumber\\
  \le D_n^{(1)}+D_{n}^{(2)} + \frac{\gk_n }2\left(
\left( S_{n}^{(1)}\right)^2 +\left( {S'}_{n}^{(2)}\right)^2+
\left( S_{n}^{(2)}\right)^2 +\left( {S'}_{n}^{(1)}\right)^2
\right).
\label{iterDn}
\end{gather} 
Since clearly $D_{n+1}$ vanishes when either of the two GW trees is
empty, one has for every $v\in[0,1]$
\begin{multline}
 \E_{n+1}^{\otimes2} \left[e^{vD_{n+1}}\right] \le
 \frac{1}{B^2} \E_{n}^{\otimes2} \left[e^{vD_{n}}
     \exp\left( \frac{v}{2}\gk_n \left( \left( S_{n}\right)^2 + \left(S'_{n}\right)^2 \right)\right)\right]^{2}
+\frac{B^2-1}{B^2} \\
      \leq \frac{1}{B^2} e^{ c_0 v \gp^n}\E_{n}^{\otimes2} 
       \left[\exp\left( v e^{c_0 v (\gp')^n} D_{n} \right)\right]^{2} +\frac{B^2-1}{B^2},
\label{iterexpDn}
\end{multline} 
where in the second inequality we expanded $e^{vD_n}$ as in Remark \ref{rem:boundexp}
and used Corollary~\ref{cor:bound2} to get the constant
$c_0>0$ for $\gp:=(2\gk)\vee \frac{4\gk}{B^2}$ and some $\varphi'\in(\gp ,1)$. 
Then we set $v_0\leq 1$ and for $n\geq 0$ define $v_{n+1}:= v_n e^{-c_0 v_n (\gp')^n} \leq v_0$. Define
$X_n:= \E_{n}^{\otimes2} \left[\exp\left( v_n D_{n} \right)\right]-1$, so that using the previous inequality
one has
\begin{equation}
 X_{n+1} \leq \frac{1 }{B^2}e^{ c_0 v_n \gp^n}(X_n+1)^2 -\frac{1}{B^2} \leq
   \frac{2 e^{ c_0v_0\gp^n}}{B^2} X_n\left( 1+\frac{X_n}{2} \right) + c v_0 \gp^n.
\end{equation}
We consider the different cases $B<B_c$, $B=B_c$ and $B>B_c$
separately, but each time we estimate from above $\E_n^{\otimes
  2}\left[ e^{v_n D_n}\right]$.  One then easily deduces Proposition
\ref{prop:bound22} using that there exists a constant $c_1$ such that
$v_n\geq c_1 v_0$, and then $\E_n^{\otimes 2}\left[ e^{c_1 v_0
    D_n}\right] \leq 1+X_n$. One concludes taking $u:=c_1 v_0$.

In the sequel we actually study the iteration
\begin{equation}
\hat X_{n+1} =  \frac{2 e^{ w_n}}{B^2} \hat X_n\left( 1+\frac{\hat X_n}{2} \right) + (c/c_0)w_n,\;\;\;\;\hat X_0=X_0
\label{iterXn}
\end{equation}
where we defined $w_n:=c_0v_0 \gp^n$. Clearly, $X_n\le \hat X_n$ for every $n$.




\smallskip
- Take $B>B_c:=\sqrt{2}$.  Let us fix some
$\Phi\in\left(\frac{2}{B^2}\vee \gp,1\right)$. One has that $X_0\leq
C_0 v_0$ and one shows easily by iteration, using \eqref{iterXn} and
the definition of $w_n$, that $\hat X_n\leq C_n \Phi^n v_0$, with
$(C_n)_{n\in\N}$ an increasing sequence satisfying
\begin{equation}
 C_{n+1} = C_n e^{w_n} \left( 1+\frac12 C_n v_0 \Phi^n \right) +c' \gp^n\Phi^{-(n+1)}
\end{equation} 
(use that $\Phi>(2/B^2)$).  Then we show that provided that $v_0$ has
been chosen small enough, $(C_n)_{n\in\N}$ is a bounded sequence.
Indeed, using that $C_n\geq C_0$ one has
\begin{multline}
C_{n+1} \leq C_n e^{w_n} \left( 1+\frac12 C_n v_0 \Phi^n + c'\Phi^{-1} C_n^{-1} (\gp/\Phi)^n\right)\\
     \leq C_n e^{w_n} \exp\left( \frac12 C_n v_0 \Phi^n \right) \exp\left( c'' (\gp/\Phi)^n \right)
      \leq A \exp\left(  \frac12 v_0 \sum_{k=0}^{n} C_k \Phi^k \right).
\end{multline}
where we noted $A:=\prod_{n=0}^{\infty} e^{w_n} e^{c'' (\gp/\Phi)^n}$, with $A<+\infty$
thanks to the definition of $w_n$ and using that $\Phi>\gp$.
It is then not difficult to see that if $v_0$ is chosen small enough,
more precisely such that $A \exp\left(  v_0 C_0 \sum_{k=0}^{n} \Phi^k \right)\leq 2C_0$, then
$C_n$ remains smaller than $2C_0$ for every $n\in\N$.
From this, one gets that $X_n\leq 2C_0 \Phi^n v_0$ for every $n$.




\medskip

- Take $B< B_c$.  The idea is that if $X_0$ is small enough,
\eqref{iterXn} can be approximated by the iteration $X_{n+1}\leq
\frac{2}{B^2} X_n$ while $X_n$ remains small.  For any fixed $n\geq
0$, one chooses $v_0= a\left( B^2/2 \right)^n$ with $a$ small (chosen
in a moment), and one has $X_0 \leq C_0 a \left( \frac{B^2}{2}
\right)^n$.  Then one shows by iteration that
\begin{equation}
\hat X_k\leq C_k a\left(  B^2/2 \right)^{n-k}
\end{equation}
for some increasing sequence
$(C_k)_{k\in\N}$ verifying
\begin{equation}
 C_{k+1} = e^{w_k} C_k \left( 1+ \frac{C_k}{2} a \left( \frac{B^2}{2} \right)^{n-k} \right) 
                   + a^{-1} \left( \frac{B^2}{2}\right)^{k+1-n} w_k.
\label{iterCn}
\end{equation} 
One then shows with the same method as in the case $B>B_c$ that $C_n$ is bounded by some constant $C$ uniformly in $n$, provided that $a$
had been chosen small enough.
Thus taking $k=n$ one has
$X_n\leq c a = c v_0 \left( 2/B^2 \right)^n$.

\medskip
- Take $B=B_c=\sqrt{2}$. The iteration \eqref{iterXn} gives
\begin{equation}
   X_{n+1} \le  e^{ w_n} X_n\left( 1+\frac{ X_n}{2} \right) + 
(c/c_0)w_n,
\label{iterXnbis}
\end{equation}
and we recall that $w_n=c_0v_0 \gp^n$.  Take $v_0=\gep^{3}$, so that
$X_{0}\leq \gep$ for $\gep$ small.  We now show that if $\gep\leq
\gep_0$ with $\gep_0$ chosen small enough, one has for all $n\leq
\frac12 \gep^{-1}$ that $X_{n}\leq \gep \left( 1+ n\gep \right)$.  We
prove this by induction. For
$n=0$ this is just because one chose $X_{0}\leq \gep$. If
 $X_n\leq \gep \left( 1+ n\gep \right)$ and $n\gep \leq 1/2$,
one has (note that $w_n\leq c_0 \gep^{3}$ for all $n$)
\begin{multline}
 X_{n+1} \leq e^{c_0  \gep^{3} } \gep  \left( 1+ n\gep \right)
                \left( 1+ \frac12 \gep \left( 1+ n\gep \right) \right) +c \gep^{3},\\
    \leq  \gep \left[  (1+c_0'\gep^{3})\left( 1+ n\gep \right) 
      \left( 1+3\gep/4 \right)  +c \gep^{2}\right] \\
    \leq \gep \left[ 1+\gep \left( n+ 3/4 + c_0'\gep^{2}+c\gep \right) \right] \leq \gep\left( 1+ (n+1)\gep \right),
\end{multline}
provided that $\gep\leq\gep_0$ with $\gep_0$ small enough. This concludes the induction step.
Thus one has that for all $n\leq \frac12 \gep^{-1}$, $X_n\leq 2 \gep $, with $\gep= v_0^{1/3}$.
\end{proof}

\section{Hierachically correlated Gaussian vectors}
\begin{lemma}
\label{lem:matrix}
Let $m(\cdot)$ be a function from $\N$ to $\R$ and for $n\in\N$ let 
Let $M:=M^{(n)}=(M_{ij})_{1\leq i,j\leq 2^{n}}$ be the $2^n \times 2^n$ matrix
with entries
$M_{ij}:=m (d(i,j))$.
Then, the eigenvectors of such
a matrix do not depend on the function $m(\cdot)$, and the eigenvalues are
\begin{equation}
\begin{array}{ll}
\gl_0 = m(0)+ \sum_{k=1}^{n} 2^{k-1} m(k)           \     ,   &  \text{ with multiplicity } 1\\
\gl_p = m(0)+ \sum_{k=1}^{n-p} 2^{k-1} m(k)
     - 2^{n-p} m({n+1-p})\ ,&  \text{ with multiplicity } 2^{p-1}, \text{ for } 1\leq p\leq n .
\end{array}
\label{eq:eigenval}
\end{equation}
\end{lemma}
This comes directly from the fact that 
\begin{equation}
 M^{(n)} =
\left(
\begin{array}{cc}
M^{(n-1)}         & \begin{array}{ccc}
                    m(n) & \cdots & m(n)\\
                    \vdots &       & \vdots \\
                    m(n) & \cdots & m(n)
                  \end{array}                     \\
\begin{array}{ccc}
  m(n) & \cdots & m(n) \\
 \vdots &       & \vdots \\
  m(n) & \cdots & m(n)
\end{array}                    &  M^{(n-1)} 

\end{array}
\right),
\end{equation}
where each block is of size $2^{n-1}$. One computes the eigenvalues:
the eigenvector $(1,\ldots,1)$ gives $\lambda_0$, the eigenvector
$(1,\ldots,1,-1,\ldots,-1)$ gives $\gl_1$. Then the eigenvectors
$(X,0)$ and $(0,X)$ with $X\neq (1,\ldots,1)$ being an eigenvector of
$M^{(n-1)}$ give all the others eigenvalues, which are the eigenvalue
associated to $X$ with $M^{(n-1)}$, but with multiplicity multiplied by
$2$.

\begin{rem}\rm
Lemma \ref{lem:matrix} shows that the spectral radius of $M^{(n)}$ is
upper bounded by 
$\sum_{p=0}^{\infty} 2^p |m(p)|$. Also,
two matrices with entries depending only on the distances
$d(i,j)$ can be codiagonalized, as the eigenvectors do not depend on the values of the entries,
and one can describe the diagonalizing orthogonal matrix $\Omega$
\begin{equation}
\label{defOmega}
 \Omega= \frac{1}{\sqrt{2^{n}}}
 \left(
\begin{array}{rrrrr}
    1     &     1    & \sqrt{2}&  0       &        \\
  \vdots  &  \vdots  &  \vdots & \vdots   & \cdots \\
    1     &     1    &-\sqrt{2}&   0      &        \\
    1     &    -1    &    0    & \sqrt{2} &        \\
  \vdots  &  \vdots  &  \vdots & \vdots   & \cdots \\
    1     &    -1    &    0    & -\sqrt{2}&        
\end{array}
\right)
\end{equation} 
such that $\gO^t K\gO = Diag\left( \gl_0,\gl_1,\gl_2,\gl_2,\ldots \right)$ with $\gl_i$ given in Lemma \ref{lem:matrix}.
\label{rem:eigenval}
\end{rem}

\medskip Let $\omega=\{\omega_i\}_{i\in\N}$ be the centered Gaussian
family with correlation structure
$\bbE[\omega_i\omega_j]=\kappa_{d(i,j)}$.  The following Proposition
gives the dependence on $\kappa_n$ of a smooth function of
$\go_1,\ldots,\go_{2^n}$:
\begin{proposition}
If  $f:\R^{2^n}\mapsto \R$ is twice differentiable and grows at most
polynomially at infinity, one has
\begin{equation}
 \frac{\partial }{\partial \gk_n} \Eo\left[f(\go_1,\ldots,\go_{2^n}) \right]  =
       \sum_{i=1}^{2^{n-1}} \sum_{j=2^{n-1}+1}^{2^n} \Eo\left[\frac{\partial^2 f}{\partial \go_i \partial \go_j}(\go) \right].
\end{equation}
\label{prop:derivEot}
\end{proposition}

\begin{proof}
Thanks to Remark \ref{rem:eigenval}, one has
\begin{equation}
 \Eo\left[f(\go_1,\ldots,\go_{2^n}) \right] = \tilde\Eo\left[ f(\gO \go)\right],
\end{equation}
with $\gO$ defined in \eqref{defOmega}, and where $\tilde\Po$ stands
for the law of a centered Gaussian vector of covariance matrix
$\gD:=Diag\left( \gl_0,\gl_1,\gl_2,\gl_2,\ldots \right)$. The
eigenvalues $\lambda_i$ and their multiplicity are given in Lemma
\ref{lem:matrix}. Then, as only $\lambda_0=\gk_0+ \sum_{k=1}^{n}
2^{k-1}\gk_k $ and $\lambda_1=\gk_0+ \sum_{k=1}^{n-1} 2^{k-1}\gk_k -
2^{n-1}\gk_n$ depend on $\gk_n$ one gets
\begin{equation}
 \frac{\partial }{\partial \gk_n} \Eo\left[f(\go) \right] =  2^{n-1} \frac{\partial }{\partial \lambda_0} \tilde\Eo\left[ f(\gO \go)\right]
    -2^{n-1} \frac{\partial }{\partial \lambda_1} \tilde\Eo\left[ f(\gO \go)\right].
\label{derivEotf1}
\end{equation} 
Then one uses the classical Gaussian fact that if $\omega$ is a 
centered Gaussian variable of variance $\sigma^2$ and  $g$ is a
differentiable function which grows at most polynomially at infinity, 
\begin{equation}
 \frac{\partial}{\partial \sigma^2}\Eo\left[  g(\go)\right] = 
       \frac12 \Eo\left[ \frac{\partial^2 g}{ \partial \go^2} (\go)\right].
\end{equation}
Plugging this result in \eqref{derivEotf1}  one gets
\begin{multline}
  \frac{1}{2^{n-1}}\frac{\partial}{ \partial \gk_n} \Eo[f(\go_1,\ldots,\go_{2^n})] \\
  =\frac12\sum_{i,j=1}^{2^{n}} \gO_{i1} \gO_{j1} \tilde\Eo\left[
\left.    \frac{\partial^2 f}{ \partial x_i \partial x_{j}}
\right|_{x=\Omega\go}
  \right]
  - \frac12\sum_{i,j=1}^{2^{n}} \gO_{i2} \gO_{j2} \tilde\Eo\left[
\left.    \frac{\partial^2 f}{ \partial x_i \partial x_{j}}
\right|_{x=\Omega\go}\right]\\
  = \frac{1}{2^{n}} \sumtwo{i,j=1}{d(i,j)=n}^{2^{n}} \Eo\left[
    \frac{\partial^2 f}{ \partial\go_i \partial\go_{j}} (\go) \right],
\end{multline}
where in the second equality we used the values of $\gO_{k1}$ and $\gO_{k2}$.
\end{proof}

\begin{rem}\rm
 With the same type of computations, since $\gO$ is explicit, one can also compute the derivative with respect to $\gk_p$ for $p\leq n$,
and after some computations, one gets
\begin{equation}
  \frac{\partial }{\partial \gk_p} \Eo\left[f(\go_1,\ldots,\go_{2^n}) \right]  =
      \frac{1}{2} \sumtwo{i,j=1}{d(i,j)=p}^{2^{n}} \Eo\left[\frac{\partial^2 f}{\partial \go_i \partial \go_j}(\go) \right].
\end{equation} 
\end{rem}

\end{appendix}

\section*{Acknowledgments}
We wish to thank Bernard Derrida for a very interesting discussion
and Giambattista Giacomin for comments on the manuscript.
This work was partly done while the authors were visiting the Department of Mathematics of the University of Roma Tre, with financial support by ERC
“Advanced Grant” PTRELSS 228032.

\end{document}